\documentclass[preprint,3p,11pt]{elsarticle}

\usepackage{amssymb}
\usepackage{amsmath}
\usepackage{graphicx}
\usepackage{booktabs} 

\usepackage{array,tabularx,tabulary,booktabs} 
\usepackage{longtable}  
\usepackage{multirow}   

\biboptions{sort&compress}

\newcommand{\E}[1]{\mathop{{\rm \bf E}\!\left\{#1\right\}}\nolimits}

\begin{document}

\begin{frontmatter}
\title{Evolving efficiency of the BRICS markets}

\author[CEMAT]{Maria V. Kulikova\corref{cor}} \ead{maria.kulikova@ist.utl.pt} \cortext[cor]{Corresponding
author.}

\author[UCT]{David R. Taylor} \ead{David.Taylor@uct.ac.za}

\author[CEMAT]{Gennady Yu. Kulikov} \ead{gennady.kulikov@tecnico.ulisboa.pt}

\address[CEMAT]{CEMAT, Instituto Superior T\'ecnico, Universidade de Lisboa, Av.~Rovisco Pais, 1049-001 Lisboa, Portugal.}
\address[UCT]{The African Institute of Financial Markets \& Risk Management, University of Cape Town, Rondebosch 7701, Cape Town, South Africa.}

\begin{abstract}
This paper investigates a time-varying version of weak-form market efficiency in the BRICS countries. A moving window test for sample autocorrelations is applied alongside a Kalman filter approach to recover the hidden dynamics of the market efficiency process through appropriate time-varying autoregressive models with both homoscedastic and heteroscedastic conditional variance. Monthly data covers the period from January 1995 to December 2020, which includes the 2008-2009 global financial crisis and the recent COVID-19 recession. The results reveal that all the BRICS stock markets were affected during both periods, but generally remained weak-form efficient, with the exception of China.
\end{abstract}


\end{frontmatter}

\section{Introduction}\label{sect1}

Many recent studies have focused on a time-varying version of weak-form market efficiency that is a key component of the so-called Adaptive Market Hypothesis (AMH), introduced by Andrew Lo in~\cite{2004:Lo}. In subsequent work~\cite{2005:Lo}, the AMH is shown to `reconcile' Eugene Fama's  efficient market hypothesis (EMH)~\cite{1965:Fama} with behavioural finance. For example, the AMH proposes that investors can make mistakes, but that they are able to learn from them and adapt their behavior accordingly. Furthermore, any measure of market efficiency is assumed to fluctuate over time where this process is influenced by market conditions and economic policy. This leads naturally to time-varying econometric models and estimation techniques, as well as the derivation of various rolling window tests for tracking the change in market efficiency, see~\cite{2005:Jagric,2014:Ghazani,2019:Boya,2020:Jiang}, and others. A comprehensive and systematic review of weak-form market efficiency literature can be found in~\cite{2011:Lim}.

One of the most common methodologies for modeling and estimating a degree of market efficiency lies in an analysis of the serial autocorrelation in an observed return series. This is frequently motivated by the following. The EMH rests on predictability, ``a market is weak-form efficient when there is no predictable profit opportunity based on the past movement of asset prices''. An immediate implication of this would be the absence of serial autocorrelation in returns, thus eliminating some degree of predictability. Consequently, autocorrelation of monthly stock returns is commonly considered as a good proxy of market inefficiency in the econometric literature when investigating the AMH in the weak sense. A simple approach to modeling any degree of efficiency over time is to assume that the return series follows an AR model with time-varying coefficients and homoscedastic conditional variance, see, for example,~\cite{2009:ItoSugiyama,2011:Kim:AMH}. However, modern econometric trends suggest the use of more sophisticated system structures to model complicated `stylized' facts such as volatility clustering via Autoregressive Conditional Heteroscedasticity (ARCH) and generalized ARCH (GARCH) models~\cite{1982:Engle,1986:Bollerslev}. Tests for weak-form  market efficiency based on time-varying AR models with various GARCH-type processes accounting for heteroscedasticity in the conditional variance have received increased attention in recent years. For instance, the test for evolving efficiency (TEE) proposed in~\cite{1997:Emerson} and formalized in~\cite{1999:Zalewska} employs a GARCH-in-Mean(1,1) process combined with an AR(1) assumption to model the conditional mean where the regression coefficients follows a random walk and should be jointly estimated when formulating a conclusion about the market's efficiency. Empirical studies based on the TEE approach include seven African stock markets in~\cite{2005:Jefferis}, eleven Arab stock markets in~\cite{2010:Abdmoulah}, as well as the Gulf Cooperation Council (the economic union consisting of all Arab states of the Persian Gulf, except Iraq) in~\cite{2016:Charfeddine}, the Russian and Czech markets in~\cite{2002:Hall,2008:Posta,2019:RJ:Kulikov}, and so on. In~\cite{2000:Rockinger,2003:Li:AMH:China}, asymmetric and threshold GARCH-type specifications are proposed in evolving market efficiency tests that allow for asymmetries in the volatility of returns in reaction to information shocks. This is referred to as the leverage effect, and is often modeled using asymmetric GARCH-type processes. For example, the empirical study in~\cite{2003:Li:AMH:China} found significant evidence of a leverage effect in the Shanghai stock market, while the Shenzhen return series show no such asymmetry.

In this paper, we investigate the transition markets of Brazil, Russia, India, China and South Africa, collectively known as the BRICS countries. Although there are numerous studies of their market efficiency, very few have tested for {\it evolving} market efficiency and no work involves modeling approaches that allow for jointly measuring the degree of time-varying efficiency of BRICS markets. In contrast to earlier contributions in the literature that explore weak-form efficiency of BRICS markets in the {\it absolute} sense, reported in~\cite{2010:Chong,2014:Mobarek,2015:Nalin,2019:AsianJournal:Kiran}, our methods and modeling techniques examine the {\it dynamics} of the degree of market efficiency and assess whether its evolution is consistent with the AMH. We apply non-parametric moving window tests for sample autocorrelations together with time-varying GARCH-type models and the Kalman filter (KF) to recover the hidden dynamics of the efficiency process. The data we have used allows us to investigate the impact of the COVID-19 pandemic on the dynamics of the efficiency of the BRICS markets. More precisely, our monthly data covers the 25 years from January 1995 to December 2020, which allows us to account for the 2008-2009 global financial crisis and the recent COVID-19 recession. We are interested in how well the BRICS markets were able to cope with the 2008-2009 global financial crisis and how they are responding to the recent economic recession provoked by the COVID-19 pandemic. We also answer the following questions: 1) Which country has managed these crises in the most stable way? And, 2) How strong was the impact of the recent COVID-19 recession on evolving weak-form efficiency of BRICS markets when compared to the crash during the 2008-2009 global financial crisis?

We may have anticipated that all BRICS stock markets would be significantly affected by the 2008-2009 global financial crisis, but all have since recovered and showed signs of weak-form efficiency by the COVID-19 pandemic. The recent global recession is an interesting case for further investigation. We foresee that China may have been severely affected by the COVID-19 pandemic as an impact of its zero-tolerance approach implemented to manage the pandemic. The Indian stock market is also predicted to suffer from the strong lockdowns at that time. In contrast, Russian and Brazilian markets seem to be less affected by the recent COVID-19 recession because of weak restrictions implemented at that period. Based on a previous study in~\cite{2021:Okorie}, which  covers the period from June 2019 to July 2020, we  may additionally anticipate that Brazilian, Russian and Indian markets exhibit a trend toward weak-form efficiency at the end of our empirical study. Other markets under investigation have not been studied in recent years. This is an interesting open problem to be covered in our work.

The rest of the paper is organized as follows. Major findings on evolving market efficiency for each country are summarized in Section~\ref{Section:survey}. A statistical description of BRICS market data is presented in Section~\ref{Section:Data:new}. Moving window tests are applied for time-varying sample autocorrelation estimation and for the analysis of changing serial dependence over time in Section~\ref{Section:Data}. By applying a Kalman filter, we investigate time-varying  autoregressive models with both homoscedastic and heteroscedastic conditional variance assumptions in order to estimate a level of evolving efficiency in Section~\ref{Section:Models}. Section~\ref{Section:conclusion} summarizes the key findings of our study.  It also outlines the problems that are still open and available for future investigation. Finally, the details of estimation algorithms and the general calibration scheme, as well as the overall setup for producing our numerical results, are summarized in the Appendices.

\section{A brief overview of the market efficiency literature for BRICS countries} \label{Section:survey}

Various forms of market efficiency in BRICS countries have already been studied. In general, an `absolute', weak-form efficiency has most often been examined. Absolute, weak-form efficiency assumes (and tests for) a constant level of efficiency and results in a simple binary conclusion. For example, the Brazilian stock market was concluded to be the most efficient among the BRIC markets in~\cite{2010:Chong}. An absolute, weak-form study based on autocorrelation and unit-root tests in~\cite{2015:Nalin} produced evidence of weak-form efficiency in Brazil and India over the period from July 1997 to December 2013. An additional cointegration test concluded that there is a long-run relationship between Indian and Chinese markets. In~\cite{2014:Mobarek}, the period from September 1995 to March 2010 was divided in three; and the study revealed a trend toward increasing, `fairly' weak-form efficiency.

The South African market is seldomly investigated together with other BRIC countries in the literature. The most recent and complete BRICS studies were in~\cite{2019:AsianJournal:Kiran,2019:Saluja}. These showed that efficiency is indeed time-varying, and serial correlated, Ljung Box and runs tests indicated that all the markets demonstrated weak-form efficiency over the period. However, the Hurst exponent results showed that only the Russian market was efficient over the entire sample period from 2000 to 2018. Finally, we note that the rejection or acceptance of the market efficiency hypothesis for each market in the BRICS group depends critically on the period for which the research has been carried out. This is in line with the time-varying, weak-form efficiency assumption of the AMH framework. Conclusions about evolving weak-form efficiency also depend on the observation frequency, i.e. on whether daily, weekly or monthly returns have been investigated. It is worth noting that the various tests developed for measuring market (in)efficiency by estimating a serial correlation typically yield a conclusion about the existence of serial correlation in daily returns. For these reasons and due to the fact that different studies deliver varied results, we decided to represent our literature review on {\it changing} market efficiency of the BRICS countries in tabular form (Table~\ref{tab:survey}). This is a select overview and only {\it evolving market efficiency} results within the AMH framework are summarized. This enables a quick comparison between findings for each country.

{\tiny
\begin{longtable}{p{0.06\textwidth}p{0.03\textwidth}p{0.09\textwidth}p{0.08\textwidth}p{0.09\textwidth}p{0.25\textwidth}p{0.25\textwidth}}
\caption{Select literature review on time-varying market efficiency tests for BRICS countries} \label{tab:survey}   \\                                                                                                                                                                                                                                                                                                                                                                                                                                                                                                                                                                                                                                                             		 \toprule
{\bf Country} & {\bf Year} & {\bf Reference} & {\bf Data and frequency} & {\bf Period}  & {\bf Test(s) for estimating evolving weak-form efficiency}  & {\bf Conclusions} \\
	\endfirsthead
\toprule
{\bf Country} & {\bf Year} & {\bf Reference} & {\bf Data and frequency} & {\bf Period}  & {\bf Test(s) for estimating evolving weak-form efficiency}  & {\bf Conclusions} \\
\toprule
	\endhead
\hline
	\multicolumn{7}{r}{Table continued\ldots} \
	\endfoot
\hline
	\endlastfoot
\toprule
Brazil & 2004 & Cajueiro, Tabak~\cite{2004:Cajueiro} & Bovespa (daily)   & Jan.1992 - Dec.2002  &  The data was filtered through an AR-GARCH procedure and the resulting Hurst exponents were calculated & $\bullet$ There is evidence of a trend towards increasing efficiency \\
& 2014 & Dourado, Tabak~\cite{2014:Dourado} & Bovespa (daily)   & Jan.1995 - Dec.2012  &  Moving subsamples with fixed size were used and an hypothesis of random walk behavior was tested & $\bullet$ The random walk hypothesis could not be rejected\\
& 2017 & Mitra {\it et al.}~\cite{2017:Mitra}  & Bovespa (monthly) & Dec.1999 - Nov.2015 &  A moving window of 24-month subsamples was used and p-values of a generalized spectral test were estimated &  $\bullet$ A departure from weak-form efficiency is detectable \\
& 2021 & Okorie, Lin~\cite{2021:Okorie}  &  Bovespa (daily) & Jun.2019 - Jul.2020 & A rolling-window, martingale difference test and a conditional heteroscedasticity test were applied &  $\bullet$ There is no evidence of a substantial change in the level of weak-form efficiency for the Brazilian stock market in the short, medium, and long terms  \\
 \cline{2-7}\\
Russia & 2000 & Rockinger, Urga~\cite{2000:Rockinger}  &  ROS (daily) & Apr.1994 - June 1999 & TEE: Time-varying AR(1) process with random walk and asymmetric GARCH coefficients was estimated&  $\bullet$ The Russian market is significantly predictable \\
& 2002 & Hall, Urga~\cite{2002:Hall} & RTSI,  ASPGEN (daily) & Sep.1995 - Mar.2000  & TEE: Time-varying AR(1) process with random walk and GARCH-in-Mean coefficients was estimated& $\bullet$ The market was initially weak-form inefficient, took approximately 2.5 years to display efficiency, and there is evidence of an increase in efficiency \\
& 2005 & Jagric {\it et al.}~\cite{2005:Jagric} & RTSI (daily) & Sep.1995 - Aug.2004 & Long-range dependence (LRD) based on the Hurst exponent over rolling subsamples was estimated & $\bullet$ The Hurst exponent is time-varying and there is evidence of an increase in efficiency \\
& 2006 & Cajueiro, Tabak~\cite{2006:Cajueiro} & RTSI (daily) & Sep.1995 - May 2004 & LRD based on the Hurst exponent and short-term predictability using the variance ratio (VR) test was estimated & $\bullet$ The VR test rejected the random walk hypothesis \\
& 2009 & Risso~\cite{2009:Risso} & RTSI (daily) & July 1997 - Dec.2007 & 20 stock markets were ranked according to an efficiency measure using a symbolic time series
analysis and Shannon entropy & $\bullet$ The Russian market is found to be
one of the most inefficient markets (ranked seventeenth) \\
& 2009 & Anatolyev~\cite{2009:Anatolyev} & RUX (weekly) & Jan.1998 - Jan.2005 & Nonparametric tests for predictability of returns were applied& $\bullet$ The hypothesis of mean predictability is rejected for the Russian market \\
& 2012 & Ivanov {\it et al.}~\cite{2012:Ivanov} & RTSI (daily) & Oct.2000 -  Aug.2010  & LRD tests and an investigation of forecasting possibilities were performed & $\bullet$ The Russian market has LRD and predictability \\
& 2014 & Bogdanova, Ivanov~\cite{2014:Bogdanova} & RTSI (daily) & Jan.2002 -  Dec.2012  &  LRD tests based on a wavelet technique and an investigation of forecasting
possibilities was performed& $\bullet$ The Russian market is characterized by a decreasing level of predictability\\
& 2019 & Kulikov {\it et al.}~\cite{2019:RJ:Kulikov}  &  RTSI (daily) & Mar.2002 - Mar.2012 & TEE: Time-varying AR(1) process with random walk and GARCH-in-Mean coefficients with nonlinear feedback was estimated& $\bullet$ The Russian market is markedly weak-form inefficient \\
& 2021 & Okorie, Lin~\cite{2021:Okorie}  & MOEX (daily) & Jun.2019 - Jul.2020 & A rolling-window, martingale difference test and a conditional heteroscedasticity test were applied &  $\bullet$ There is evidence of an increase in efficiency after the COVID-19 outbreak  \\
 \cline{2-7}\\
India  & 2014 & Hiremath, Kumari~\cite{2014:Hiremath}  &  Sensex, Nifty (daily) & Jan.1991 - Mar.2013 & AMH: Linear and nonlinear methods were applied to empirically test the hypothesis  & $\bullet$ There is evidence of an increase in efficiency \\
& 2015 & Mishra {\it et al.}~\cite{2015:Mishra}  & Six main indices (monthly) & Jan.1995 - Dec.2013 & Three unit-root tests with two structural breaks, accounting for GARCH effects were applied & $\bullet$ The Indian indices are mean-reverting \\
& 2017 & Mitra {\it et al.}~\cite{2017:Mitra}  & BSESN (monthly) & Dec.1999 - Nov.2015 &  A moving window of 24-month subsamples was used and p-values of a generalized spectral test were estimated &  $\bullet$ A departure from weak-form efficiency is detectable \\
& 2020 & Bhuyan {\it et al.}~\cite{2020:Bhuyan}  & Sensex, Nifty (daily) & Apr.1979 - Oct.2019 & The automatic portmanteau ratio and wild bootstrap automatic variance ratio test statistics were estimated in a rolling window framework &  $\bullet$
 The BSE and NSE are informationally inefficient in the weak-form sense and the degree of predictability evolves over the period \\
& 2021 & Okorie, Lin~\cite{2021:Okorie}  & BSE Sensex 30 (daily) & Jun.2019 - Jul.2020 & A rolling-window, martingale difference test and a conditional heteroscedasticity test were applied &  $\bullet$ There is evidence of an increase in efficiency after the COVID-19 outbreak  \\
 \cline{2-7}\\
China & 2003 & Xiao-Ming Li~\cite{2003:Li:AMH:China} & Shanghai and Shenzhen SE Composite (daily) & Apr.1991 -  Jan.2001 & TEE: Time-varying AR(n), $n=1,2,3$, processes with random walk and asymmetric GARCH coefficients were estimated&  $\bullet$ A leverage effect is detected in the Shanghai market, but not in the Shenzhen market \\
& 2003 & Chen, Hong~\cite{2003:Chen} & Shanghai and Shenzhen major indices (daily) & Dec.1990 -  Oct.2002 & A generalized spectral derivative test in the presence of volatility clustering was applied &  $\bullet$  The Shanghai and Shenzhen markets are significantly weak-form inefficient, although there is evidence of a trend towards increasing efficiency\\
& 2017 & Shi {\it et al.}~\cite{2017:Shi}  & SHSE, SZSE (daily, monthly) & Dec.1990 - Sep.2015 & The wild bootstrap automatic variance ratio test and the generalized spectral test were applied &  $\bullet$ The predictability of returns varies over time and significant predictability is observed during market turmoil \\
& 2017 & Mitra {\it et al.}~\cite{2017:Mitra}  & BSESN (monthly) & Dec.1999 - Nov.2015 &  A moving window of 24-month subsamples was used and p-values of a generalized spectral test were estimated &  $\bullet$ A departure from weak-form efficiency is detectable \\
& 2019 & Xiong {\it et al.}~\cite{2019:Xiong}  & SSE50, SSE180, CSI300, ChiNext (daily) & inception - Dec.2015 &  Subsample and rolling window analyses with GARCH models were applied &  $\bullet$ The results suggest that the AMH is a better explanation of the Chinese market \\
& 2020 & Jiang, Li~\cite{2020:Jiang} & Shanghai Composite Index (weekly) & Jan.1997 -  Apr.2019 & Evolving market efficiency was estimated using quantile autoregression &  $\bullet$ The degree of market efficiency varies dramatically as the quantile level changes. The results are largely due to the immature financial market environment in China \\
 \cline{2-7}\\
South Africa  & 2004 & Jefferis, Smith~\cite{2004:Jefferis} & TOP40 {\it et al} (weekly)  & Jan.1993 - Mar.2001  & VR test and TEE: Time-varying AR(1) process with random walk and GARCH-in-Mean coefficients was estimated & $\bullet$ The TEE suggests that the TOP40 index is weak-form efficient. The VR test does not reject the random walk hypothesis for the TOP40.\\
 & 2005 & Jefferis, Smith~\cite{2005:Jefferis} & JSE ALSI (weekly) & Jan.1990 -  June 2001  & TEE: Time-varying AR(1) process with random walk and GARCH-in-Mean coefficients was estimated  & $\bullet$ The JSE exhibits fairly constant weak-form efficiency and shows no sign of changing\\
  & 2009 & Morris {\it et al.}~\cite{2009:Morris} &  TOP40 (daily)  & Jan.2005 -  Dec.2007  & Autoregressive fractionally integrated moving average was estimated and wavelet analysis was performed to test for the existence of long term memory & $\bullet$  Share prices have long term memory and price changes are correlated over time. There is evidence that the JSE is not efficient, even in its weak-form \\
  & 2012 & Kruger  {\it et al.}~\cite{2012:Kruger} & 109 shares from JSE (daily)  & Jan.2002 - Dec.2009  & A number of moving window tests were performed to test for nonlinear serial dependence & $\bullet$ Evidence of linear and nonlinear dependence is found. The JSE is (for the most part) an efficient market with brief periods of inefficiency \\
  & 2012 & Bonga-Bonga~\cite{2012:Bonga} &  JSE indices (weekly)  & Mar.1995 - Dec.2009  & TEE: Time-varying AR(1) process with random walk and GARCH coefficients was estimated & $\bullet$ The JSE is weak-form efficient \\
  & 2016 & Noakes, Rajaratnam~\cite{2016:Noakes} & TOP40, JSE shares, JSE ALSI (daily)  & Mar.2005 - Dec.2009  & An overlapping, serial correlation test was adapted to include adjustments for thin trading  & $\bullet$  There is evidence that many shares on the JSE were inefficient during the 2008-2009 global financial crisis, when compared to the stable period\\
  & 2018 & Heymans, Santana~\cite{2018:Heymans} & JSE indices (daily)  & July 1997 - Mar.2015  & Rolling window VR tests were performed & $\bullet$ The JSE all share index is weak-form efficient \\
  & 2019 & Kulikov {\it et al.}~\cite{2019:RJ:Kulikov}  &  TOP40 (daily) & Mar.2002 - Mar.2012 & TEE: Time-varying AR(1) with random walk and GARCH-in-Mean coefficients with nonlinear feedback was estimated &  $\bullet$ The JSE TOP40 index is weak-form efficient \\
\bottomrule
\end{longtable}
}

\section{Statistical description of BRICS markets data} \label{Section:Data:new}

We consider data from the five BRICS markets. The data we examine are the benchmark BOVESPA index of about 70 stocks traded on the Brazilian Stock Market, the Russia Trading System Index (RTSI) containing 50 Russian stocks traded on the Moscow Exchange, the BSE SENSEX, which is also known as the S\&P Bombay Stock Exchange Sensitive Index, consisting of 30 well-established and financially sound companies listed on the Bombay Stock Exchange, the Shanghai SE Composite Price Index, and the JSE TOP40, which comprises the 40 largest and most liquid stocks listed on the Johannesburg Stock Exchange.

Model estimation is based on monthly closing prices from 1995 to 2020. We use log-returns calculated in the usual way on a continuously compounded basis, i.e. $y_t = \ln S_t - \ln S_{t-1}$ where $S_t$ is the closing price at time $t$, $t=1,\ldots, N$. We note that the examination of daily log-returns often yields a conclusion of weak-form inefficiency because the tests examine serial autocorrelation, which is prevalent in daily returns. For this reason, monthly observations are preferable for an appropriate study of evolving market efficiency; see also the discussion in~\cite{2004:Lo,2009:ItoSugiyama}.

\begin{table}[ht!]
{\scriptsize
\caption{Summary statistics of monthly log returns for BRICS countries} \label{tab:stats}                                                                                                                                                                                                                                                                                                                                                                                                                                                                                                                                                                                                                                                                \begin{tabular}{lccccc}
\toprule
{\bf Country} & {\bf Brazil} & {\bf Russia} &  {\bf India} & {\bf China} & {\bf S. Africa}\\
\toprule
Index       & BOVESPA & RTSI &  BSE SENSEX & Shanghai SE & TOP40 \\
Period      & Jan.1995 - Dec.2020 & Sep.1995 - Dec.2020 & Jan.1995 - Dec.2020 & Jan.1995 - Dec.2020 & Jun.1995 - Dec.2020  \\
Sample size, $N$ & 311 & 303 & 311 & 311 & 306 \\
\hline
Mean, $\hat \mu$   & 0.01100 &  0.00917 &  0.00830 &  0.00585 &  0.00802 \\
Median             & 0.01305 &  0.01328 &  0.01030 &  0.00632 &  0.01038 \\
Std, $\hat \sigma$ & 0.08648 &  0.13147  &  0.06870 & 0.07844 &  0.05572 \\
Skewness &  -1.15623 & -1.10107 &  -0.47748 &  -0.14555 &  -0.88107 \\
Kurtosis &   7.85301 &   9.20392 & 4.44516 & 4.65687 &  7.57658 \\
\hline
$\hat \rho_1$    &  0.0204 &  0.2015 & 0.0150 & 0.0591 &   -0.0594 \\
$\hat \rho_{5}$  & -0.1118 & -0.0217 & -0.0289 & 0.0526 &  -0.1347 \\
$\hat \rho_{10}$ &  0.0896 & -0.0107 &  0.0476 &  -0.0526 & -0.0311 \\
$\hat \rho_{15}$ & -0.1449 & -0.1219 &  0.0149 &  -0.0213 & -0.0223 \\
\hline
$Q(1)$    & 0.1311 & 12.4234 &  0.0702 & 1.0969  &  1.0920 \\
$p$-value & 0.7173 &  0.0004 & 0.7911 &  0.2949   & 0.2960 \\[5pt]
$Q(5)$    &  6.3143 &  16.7297 & 2.1842 &  14.3708 &  10.6288 \\
$p$-value &  0.2768 &   0.0050 & 0.8231  & 0.0134 &  0.0593 \\[5pt]
$Q(10)$   &  11.4713  &  21.5504  &  5.0655 & 18.3408 &  15.5251 \\
$p$-value &  0.3220  &  0.0176  &  0.8868 & 0.0495 &  0.1141 \\[5pt]
$Q(15)$   &  20.8473  & 28.5289  &  8.2825 & 31.3590 &  15.9709  \\
$p$-value &  0.1418  &  0.0185 &  0.9120   &   0.0079   & 0.3840 \\
\hline
Ljung-Box Q-tests & $H_0$ is not  & $H_0$ is rejected & $H_0$ is not & $H_0$ is rejected & $H_0$ is not \\
up to lag 20      & rejected      & for some $l$      & rejected     & for some $l$      & rejected  \\
\bottomrule
\end{tabular}}
\end{table}

The summary statistics in Table~\ref{tab:stats} lead to the following conclusions. It is clear that the Russian market is the most volatile. The negative skewness for all five series indicates that their density functions either have a longer or fatter left tail than right tail. The largest skewness values are observed for the Brazilian, Russian and South African markets. For these three samples, the kurtosis values greatly exceed that of the normal distribution. This is a clear indication of a leptokurtic distribution. The sample mean is slightly less than the median for the Brazilian and Chinese series, and substantially less for the Russian, Indian and South African. This confirms a left-skewed shape for their distributions and is the result of the mean being pulled down by the long left tail.

Table~\ref{tab:stats} also summarizes the sample autocorrelation coefficients  $\hat \rho_l$ at lag $l$, $l=1,5,10,15$, as well as the corresponding Ljung-Box statistics $Q(l)$ from a Q-test for residual autocorrelation. This is a joint test for the hypothesis that the first $l$ autocorrelation coefficients are equal to zero. We test all five samples up to lag $20$ and conclude that for the Brazilian, Indian and South African markets the null hypothesis that the residuals are uncorrelated is not rejected at either a 1\% or 5\% significance level. This gives an initial indication that the Brazilian, Indian and South African markets are weak-form efficient in the absolute sense over the entire time period examined. In contrast, the Russian and Chinese markets cannot be considered as absolutely weak-form efficient even within this preliminary study since the results of the Ljung-Box Q-test for residual autocorrelation suggest rejecting the null hypothesis, i.e. there exists evidence of serial autocorrelation in the monthly returns for the Russian and Chinese series. In the next section, we explore the {\it evolving} efficiency of BRICS markets, i.e. when the level of market efficiency is assumed to be time-varying.

\section{Rolling window tests for evolving efficiency of BRICS markets} \label{Section:Data}

By definition, ``a market is weak-form efficient when there is no predictable profit opportunity based on the past movement of asset prices''. Any sign of serial autocorrelation in an observed return series contradicts this definition and denies weak-form efficiency in the market under examination. To study {\it evolving} market efficiency, various rolling window tests are traditionally developed in econometric literature. In this paper, we follow the approach suggested in~\cite{2004:Lo,2009:ItoSugiyama} to examine the evolving efficiency of BRICS markets. More precisely, we apply the moving (rolling) window method for calculating the sample autocorrelation coefficient at lag $l=1$ as an initial sign of market inefficiency. Let $w$ denote the window size, where $w<N$ and $N$ is the number of observations. Define $N-w+1$ sub-samples of size $w$ by $\{y_{t-w+1}, \ldots, y_t \}$ for $t = w, \ldots, N$ and then compute the first order autocorrelation for each sub sample. The window size is set to $w = 80$ and the corresponding $1$\% confidence bounds are calculated. The market is weak-form efficient when $\hat \rho_1$ equals zero (within the chosen confidence interval). To enhance any insight, we also provide the rolling window Ljung-Box Q-test results for residual autocorrelation at lag $l=1$ and, concurrently, plot the time-varying $p$-values. For $p$-values less than $\alpha$, say $\alpha = 0.01$ or $\alpha = 0.05$, the null hypothesis of uncorrelated residuals may be rejected at the related significance level $\alpha$. This implies some degree of inefficiency exists during that time period.

\begin{figure}[th!]
\begin{tabular}{cc}
\includegraphics[width = 0.5\textwidth]{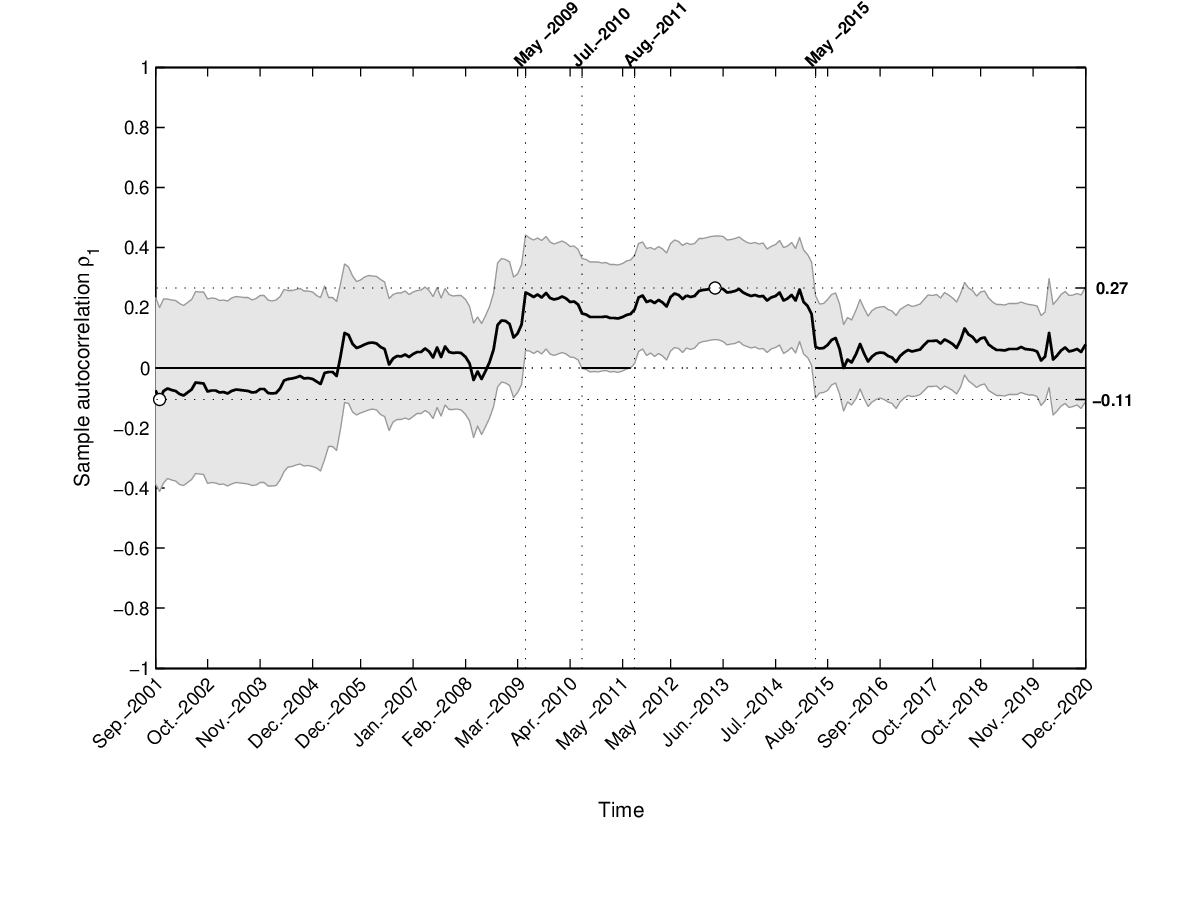} &
\includegraphics[width = 0.5\textwidth]{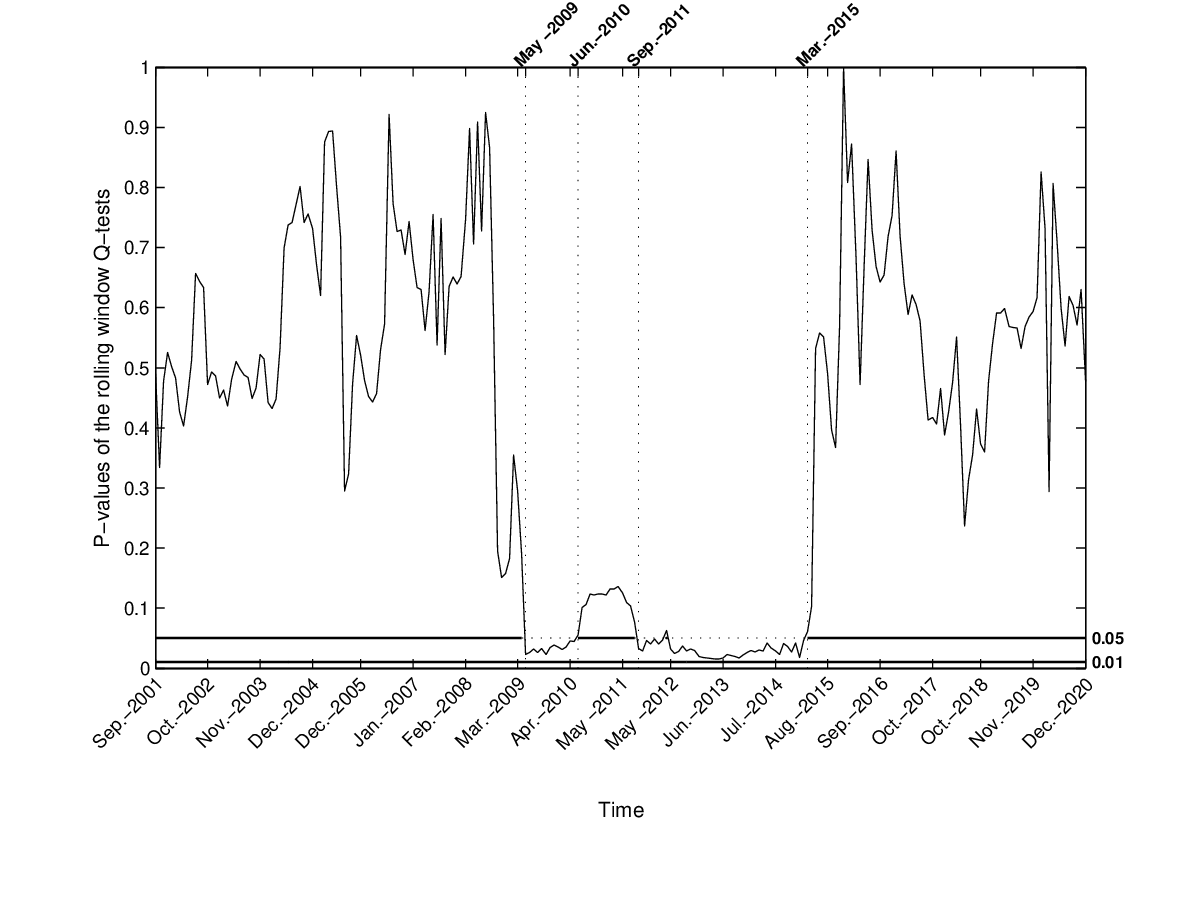}
\end{tabular}
\caption{The sample autocorrelation coefficient $\hat \rho_1$ computed by the moving window method for the Brazilian stock market with $1$\% confidence bounds (left) and the $p$-values of the rolling window Q-test for residual autocorrelation (right).}
\label{fig:MW:Brazil}
\end{figure}

We begin our discussion with an empirical study of the Brazilian market, the results of which are illustrated in Fig.~\ref{fig:MW:Brazil}. The left-hand graph represents the sample autocorrelation coefficient $\hat \rho_1$ computed by the moving window method, which is a proxy for an evolving market efficiency level. The associated $1$\% confidence bounds are also illustrated. Strong evidence of the time-varying nature of weak-form market efficiency is clearly seen from the graph. The estimated level of market efficiency fluctuates around the critical level of zero except during the period related to the 2008-2009 global financial crisis. A window size of $w=80$ months allows us to monitor the process from September 2001 to December 2020. During this period, the highest level of inefficiency is observed on April 2013 with an estimated value of $\hat \rho_1 =  0.2656$. From the estimated dynamics, we see that the process breached the efficiency level in May 2009 and remained inefficient for nearly six years. In mid-2010, the market briefly approached weak-form efficiency, but a return to weak-form efficiency only occurred in 2015.

The right-hand graph of Fig.~\ref{fig:MW:Brazil} illustrates the $p$-values of the related rolling window Ljung-Box Q-test for residual autocorrelation at lag $l=1$. This gives us an alternative version of the market efficiency process. It is clear that the Brazilian market is always weak-form efficient at a 1\% significance level, but two periods of inefficiency are detected at the 5\% level. These are the periods from May 2009 to June 2010 and from September 2011 to March 2015.

Finally, it is interesting to note that the most recent recession related to the COVID-19  pandemic resulted in no significant change in the estimated market efficiency level. From the right-hand graph of Fig.~\ref{fig:MW:Brazil}, the market is weak-form efficient from May 2015 to December 2020 with the estimated level of efficiency remaining close to the critical level of zero in the corresponding period of the left-hand graph. It is also worth noting that our conclusion concerning the COVID-19 recession concurs with the recently published results in~\cite{2021:Okorie}, where different methods were used to investigate the AMH in the Brazilian market from June 2019 to July 2020, and no substantial change in the level of efficiency was found.

\begin{figure}[th!]
\begin{tabular}{cc}
\includegraphics[width = 0.5\textwidth]{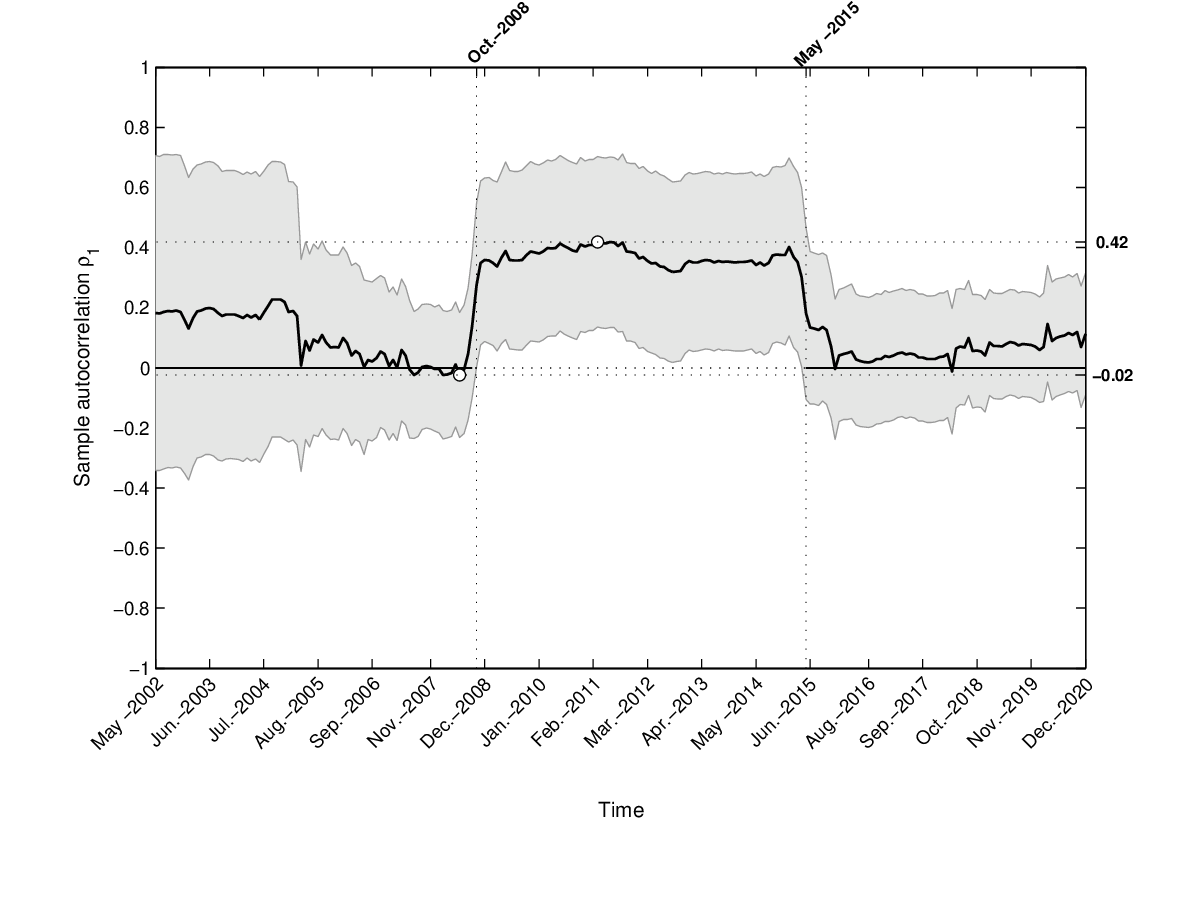} &
\includegraphics[width = 0.5\textwidth]{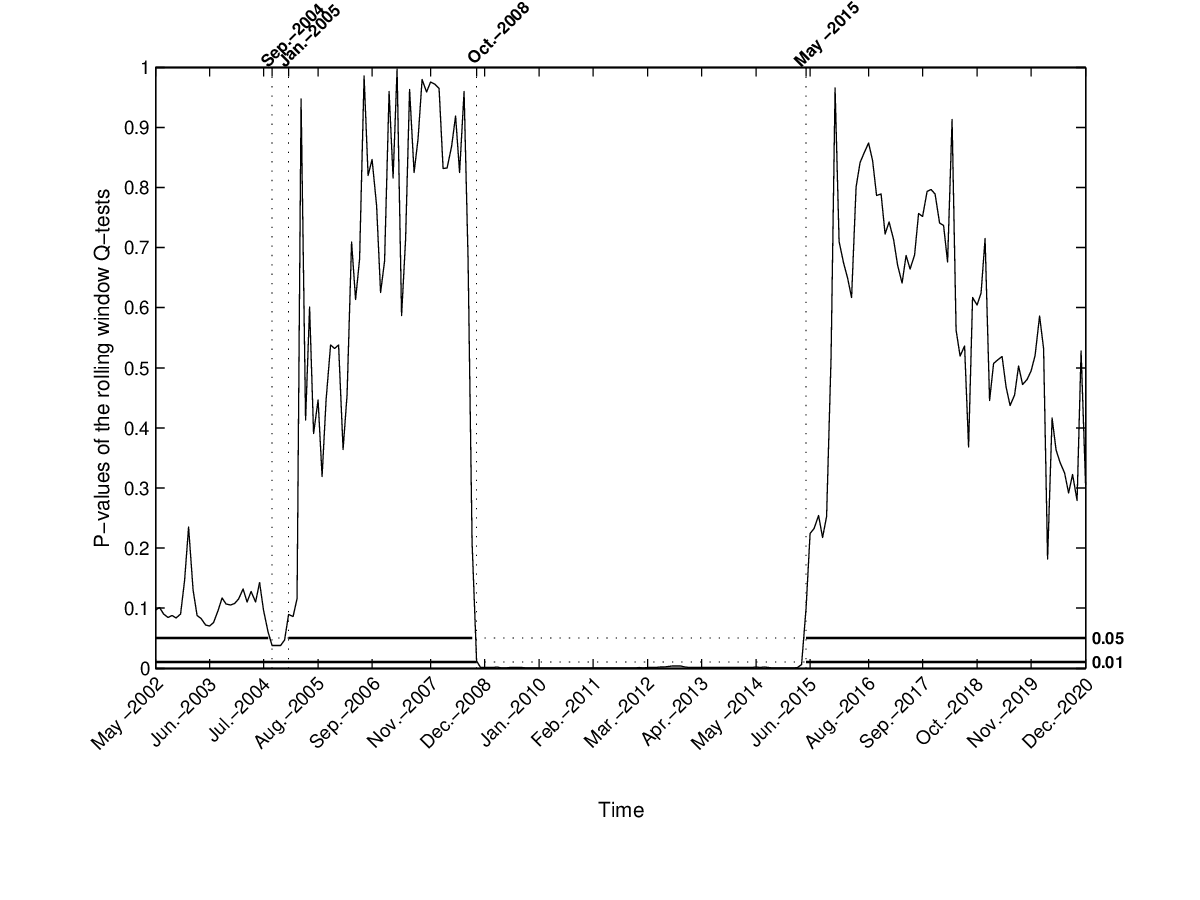}
\end{tabular}
\caption{The sample autocorrelation coefficient $\hat \rho_1$ computed by the moving window method for the Russian stock market with $1$\% confidence bounds (left) and the $p$-values of the rolling window Q-test for residual autocorrelation (right).}
\label{fig:MW:Russia}
\end{figure}

A similar pattern in the related efficiency level dynamics is observed for the Russian market, the results of which are illustrated in Fig.~\ref{fig:MW:Russia}. More precisely, we clearly see a large increase in the degree of inefficiency in the period related to the 2008-2009 global financial crisis and a long recovery process thereafter. In slight contrast to the Brazilian case, the Russian market is inefficient at both 1\% and 5\% significance levels for a period from October 2008 to May 2015. At the 5\% level, we also observe a short period of inefficiency from September 2004 to January 2005.

The maximum estimated level of inefficiency is higher than for the Brazilian index: the sample autocorrelation coefficient is estimated to be $\hat \rho_1 = 0.4190$ in March 2011. In fact, this is the highest level of inefficiency among all the countries we examined. The Russian market also has the longest recovery period from the 2008-2009 global financial crisis. Russia showed a slight trend towards weak-form efficiency in 2012-2013, but was not able to recover entirely, and only returned to efficiency at 1\% and 5\% significance levels in June 2015.

Thereafter, the estimated sample autocorrelation coefficient remains close to the critical level of zero. Similarly to the Brazilian market, no significant impact of the COVID-19 recession is observable. Both the Brazilian and Russian markets show weak-form efficiency from mid-2015 to the end of our study in December 2020. Our conclusion about the impact of COVID-19 on the Russian market again coincides with the recently published results in~\cite{2021:Okorie} obtained by alternative AMH methods.

\begin{figure}[th!]
\begin{tabular}{cc}
\includegraphics[width = 0.5\textwidth]{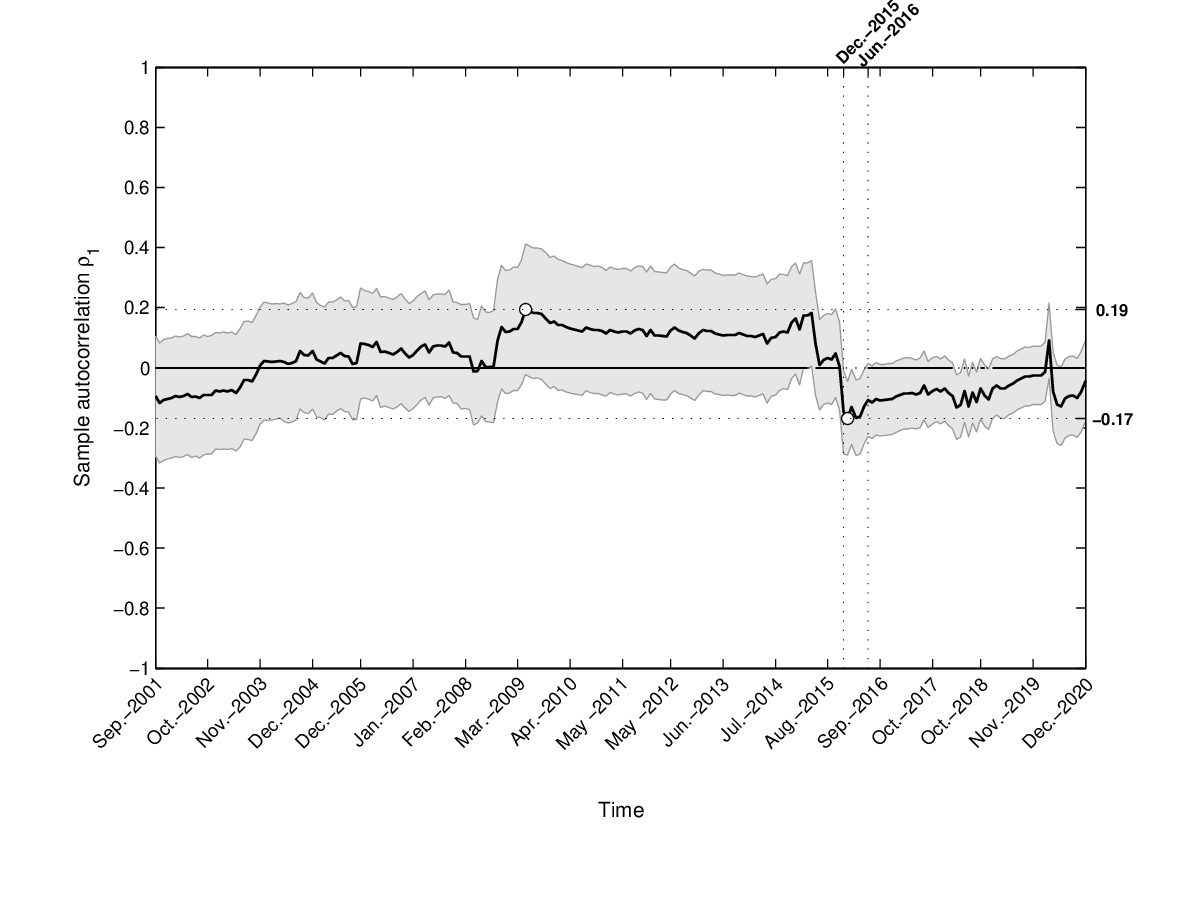} &
\includegraphics[width = 0.5\textwidth]{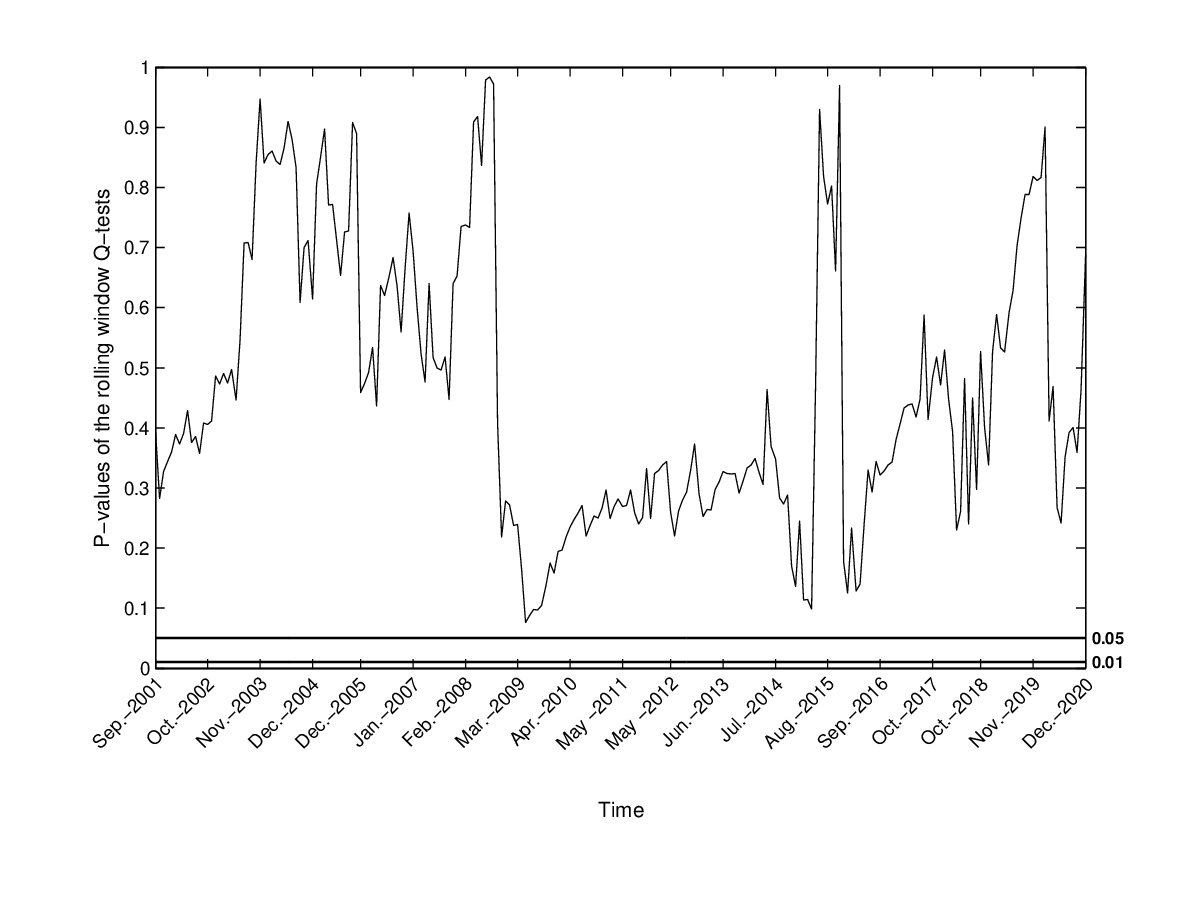}
\end{tabular}
\caption{The sample autocorrelation coefficient $\hat \rho_1$ computed by the moving window method for the Indian stock market with $1$\% confidence bounds (left) and the $p$-values of the rolling window Q-test for residual autocorrelation (right).}
\label{fig:MW:India}
\end{figure}

Next, the results of the Indian empirical study are illustrated in Fig.~\ref{fig:MW:India}. In contrast to the Brazilian and Russian cases discussed above, the Indian market behaved well during and after the 2008-2009 global financial crisis. It remains weak-form efficient for the entire period under study from September 2001 to the end of 2020 at both 1\% and 5\% significance levels, except for a short interval at the end of 2015. The estimated sample autocorrelation coefficient in the left-hand graph of Fig.~\ref{fig:MW:India} confirms a period of turbulence from December 2015 till June 2016. The maximum degree of inefficiency in the period is estimated to be $\hat \rho_1 =  0.1946$ in May 2009. This is the lowest level of inefficiency among all the countries we examined.

It is interesting to note that the recent COVID-19 pandemic had a greater effect on the Indian market than the Brazilian and Russian cases discussed above. We clearly observe an abrupt movement toward inefficiency at the end of 2019. In addition, the proxy for a time-varying level of market efficiency illustrated by the left graph  in Fig.~\ref{fig:MW:India} is highly volatile in comparison to the Brazilian and Russian cases. However, the concurrent, rolling window Q-test for residual autocorrelation does not reject the null hypothesis at either the 1\% or the 5\% significance level, i.e. the Indian market is still weak-form efficient at reasonable significance levels during the 2019-2020 COVID-19 recession, but highly unstable. Moreover, India shows a strong trend towards increased weak-form efficiency during the last year of our study. Again, this concurs with the results in \cite{2021:Okorie}.

\begin{figure}[th!]
\begin{tabular}{cc}
\includegraphics[width = 0.5\textwidth]{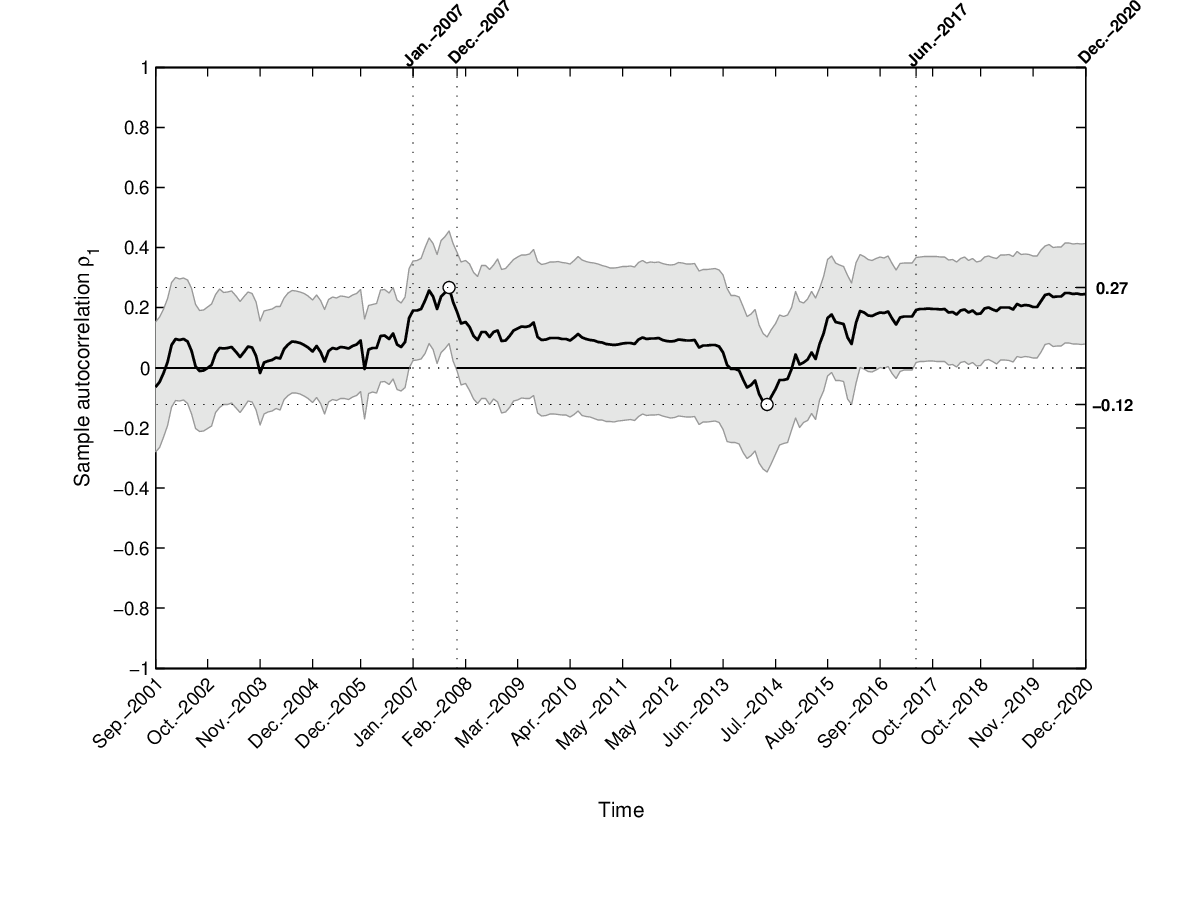} &
\includegraphics[width = 0.5\textwidth]{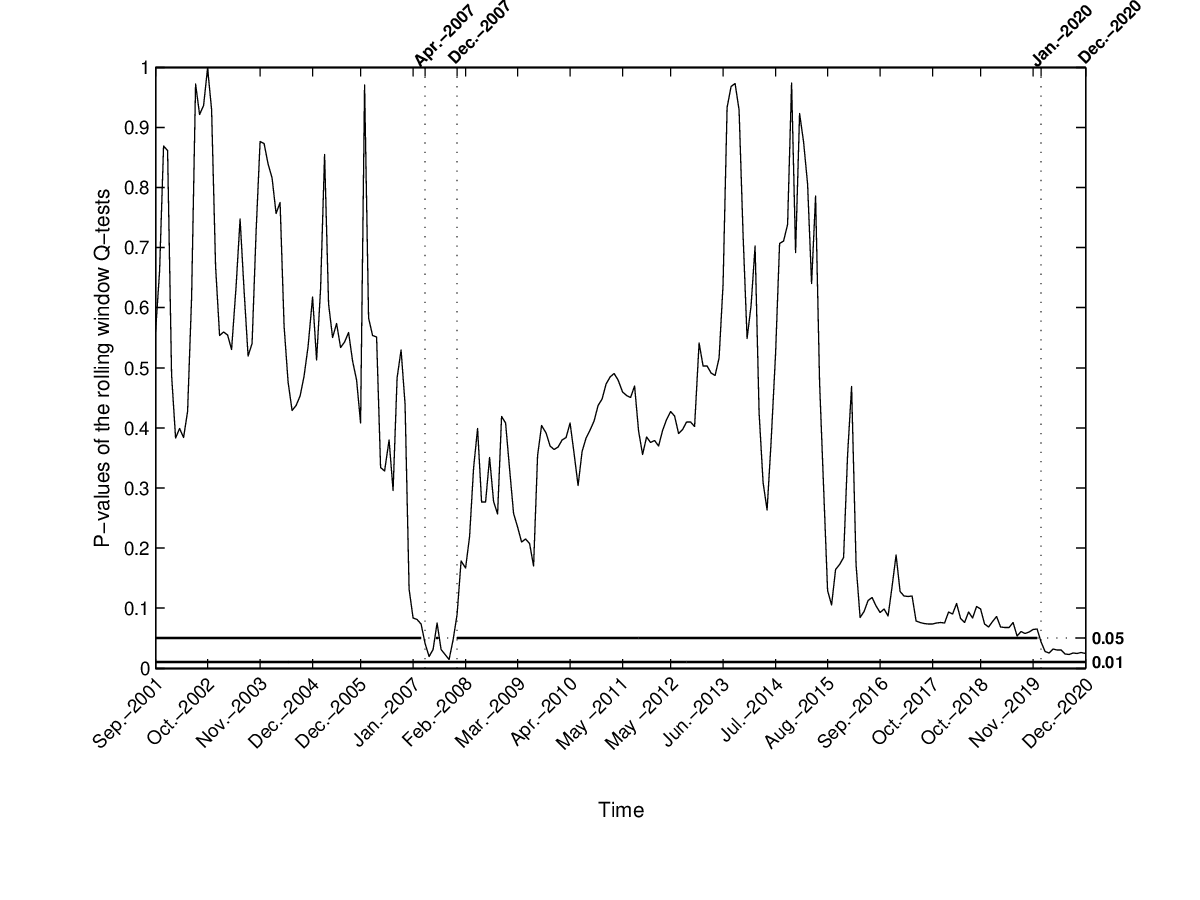}
\end{tabular}
\caption{The sample autocorrelation coefficient $\hat \rho_1$ computed by the moving window method for the Chinese stock market with $1$\% confidence bounds (left) and the $p$-values of the rolling window Q-test for residual autocorrelation (right).}
\label{fig:MW:China}
\end{figure}

Chinese markets have been most affected by the COVID-19 pandemic, where a dramatic impact on market efficiency is observable at the end of 2019. From the right-hand graph of Fig.~\ref{fig:MW:China}, we see that the Chinese market is inefficient at a 5\% significance level by January 2020 and is trending further away. The left-hand plot at Fig.~\ref{fig:MW:China} displays the dynamics of the evolving efficiency process reflected in the sample autocorrelation coefficient. As can be seen, Chinese markets have become gradually less efficient since 2014; a process that has been exacerbated by the COVID-19 pandemic. The  estimated degree of inefficiency at the end of our study $\hat \rho_1 = 0.2700$ is very close to the peak value during the 2008-2009 global financial crisis. Most importantly, China displays an opposing trend away from weak-form inefficiency in comparison to the other BRICS markets, despite negotiating the 2008-2009 global financial crisis reasonably well.

\begin{figure}[th!]
\begin{tabular}{cc}
\includegraphics[width = 0.5\textwidth]{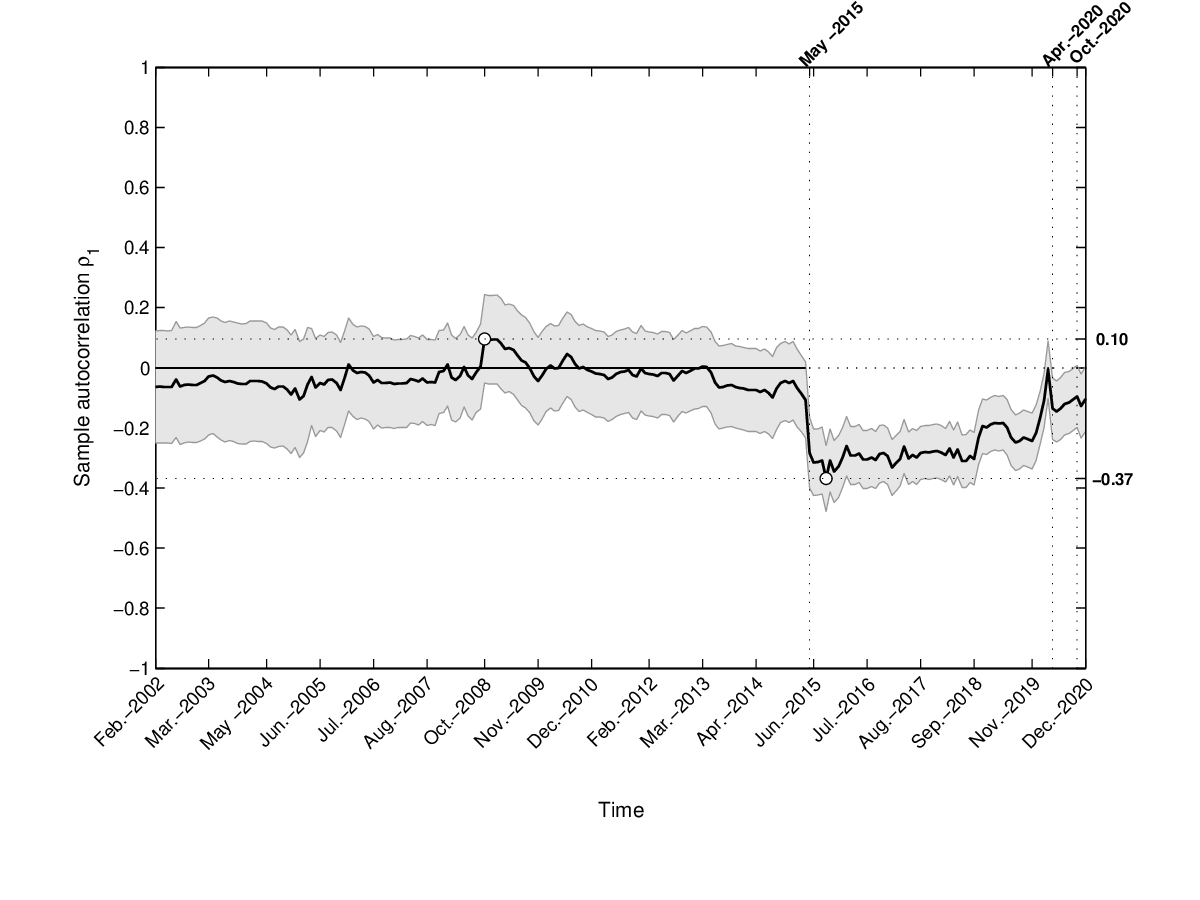} &
\includegraphics[width = 0.5\textwidth]{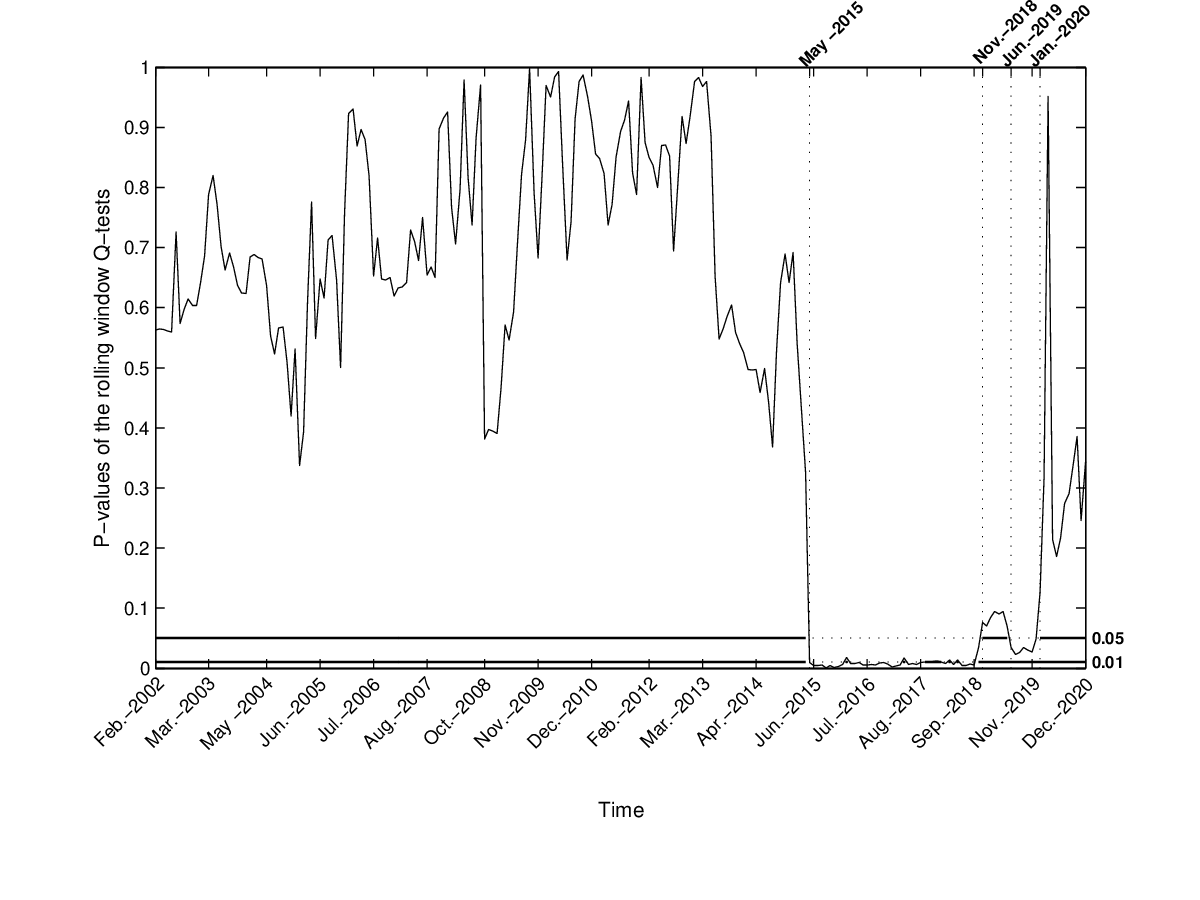}
\end{tabular}
\caption{The sample autocorrelation coefficient $\hat \rho_1$ computed by the moving window method for the South African stock market with $1$\% confidence bounds (left) and the $p$-values of the rolling window Q-test for residual autocorrelation (right).}
\label{fig:MW:SA}
\end{figure}

The South African case study illustrated in Fig.~\ref{fig:MW:SA} is interesting for the following reasons. Given the estimated dynamics in the left-hand graph, we conclude that South Africa was the most efficient and stable country in the BRICS group from February 2002 until 2015, when the situation drastically transformed. The estimated degree of market efficiency was approximately zero throughout, with a maximum value of $\hat \rho_1 =  0.10$ in December 2008. This is the smallest deviation from the critical level (efficiency level) of zero among all BRICS countries between 2002 and 2015.

However, the political crisis in South Africa in 2015 changed the situation fundamentally. The maximum degree of estimated inefficiency occurs in September 2015, with a value of $\hat \rho_1 = -0.3684$. The right-hand graph indicates an inefficiency regime at both 1\% and 5\% significance levels from July 2015 to November 2018, and at a 5\% significance level from June 2019 to the end of 2019. The large fluctuation observed in 2019 seems unrelated to the COVID-19 pandemic. Hence, we may conclude that, in contrast to the Chinese case study, South Africa has not been overly affected by the COVID-19 pandemic. In fact, the South African case is very similar to the Indian case in that an abrupt movement towards inefficiency is observed in April 2020, but the rolling window Q-test for residual autocorrelation does not reject the null hypothesis at either the 1\% or the 5\% significance levels. So, South Africa may still be considered to be a weak-form efficient in 2020, and, in fact, trends towards weak-form efficiency during 2020.

In summary, we make the following conclusions. The major findings reveal that all BRICS stock markets exhibit a trend toward increased weak-form efficiency, except China. Chinese markets display the opposite trend and have gradually become less efficient since 2017. China has also been severely affected by the COVID-19 recession. Indian and South African markets experienced significant turbulence in 2019-2020, but they show strong recovery towards weak-form efficiency in the last year. Brazilian and Russian markets have been less concerned by the COVID-19 recession, but were greatly affected by the 2008-2009 global financial crisis. Their recovery has been lingering, while China, India and South Africa were only briefly affected and were able to recover quickly. Overall, the Russian market is the most volatile and least stable, although it shows a trend toward increased weak-form efficiency in the past few years. India displays the most stable evolving efficiency process with the smallest degree of inefficiency.

\section{Modeling techniques for evolving efficiency estimation and application to BRICS} \label{Section:Models}

Modeling techniques are widely used for estimating the hidden market efficiency process as well as for designing the forecasting methods. Recall, an autocorrelation of monthly stock returns is commonly considered as a good proxy of market inefficiency in econometric literature when we deal with the AMH in the weak sense. Thus, a common approach to modeling and tracking the (time-varying) level of market efficiency within the AMH framework is to apply an autoregressive (AR) model with time-dependent coefficients to a chosen history of returns; see, for example,~\cite{1997:Emerson,1999:Zalewska,2000:Rockinger,2002:Hall,2003:Li:AMH:China,2008:Posta,2019:RJ:Kulikov} and others.
 Here, we discuss this methodology in details.

An AR process of order $n$, i.e. AR($n$), is conventionally written as:
\begin{equation}
y_t = \sum \limits_{i=1}^{n} \beta_i y_{t-i} + \varepsilon_t, \quad \varepsilon_t \sim IID(0,\sigma_{\varepsilon}^2)  \label{AR:n}
\end{equation}
where $y_t$ is the log return series and $\beta_i$, $i=1,\ldots,n$ are (constant) parameters of the model.

Under the EMH, tests for weak-form market efficiency based on an AR($n$) process require the following condition(s) to hold: $\beta_1=\beta_2= \ldots =\beta_n=0$. If this condition is met, equation~\eqref{AR:n} reduces to $y_t = \varepsilon_t$, which implies that the returns are independently distributed and that, consequently, the historical price information cannot provide consistent profit opportunities. This aligns with the definition of weak-form market efficiency proposed in~\cite{1965:Fama}. In reality, the assumption $\varepsilon_t \sim IID(0,\sigma_{\varepsilon}^2)$ is too strong for testing weak-form market efficiency. Following the discussion in~\cite{2011:Lim}, it is sufficient to assume that $\varepsilon_t$ is a white noise process in~\eqref{AR:n}. If the data suggests that $\beta_1=\beta_2= \ldots =\beta_n=0$, then this implies that the returns are serially uncorrelated in time, and again that the historical price information cannot provide profit opportunities, i.e. the market is weak-form efficient by definition. Finally, model~\eqref{AR:n} may include a constant trend $\beta_0$, but the test for market efficiency remains $\beta_1=\beta_2= \ldots =\beta_n=0$. This condition then yields $y_t = \beta_0 + \varepsilon_t$ where $\varepsilon_t$ is a white noise process, i.e. the log return is a white noise process about the mean and, hence, the returns are still serially uncorrelated in time. In this paper, we analyze mean-adjusted returns and, hence, the coefficient $\beta_0$ is omitted in the models presented below.

In contrast to the EMH, the AMH concept implies a {\it time-varying} AR model. Here we explore only AR models in the presence of Gaussian uncertainties, i.e.
\begin{align}
y_t & = \sum \limits_{i=1}^{n} \beta_{i,t} y_{t-i} + \varepsilon_t, & \varepsilon_t & \sim {\cal N}(0,\sigma_{\varepsilon}^2) \label{AR:n:time}
\end{align}
where $\beta_{i,t}:=\beta_i(t_k)$ and $\varepsilon_t$ is a Gaussian white noise process with zero mean and variance $\sigma_{\varepsilon}^2>0$. The dynamics of market efficiency is defined by the time-dependent regression coefficients $\beta_{i,t}$, $i=1, \ldots, n$, which should be modeled and estimated.

A random walk specification in \eqref{AR:n:time} for modeling the dynamics of $\beta_{i,t}$, $i=1,\ldots, n$ is the traditional model employed in the econometric literature dedicated to evolving market efficiency,
\begin{align}
\beta_{i,t} & = \beta_{i,t-1} + w_{i,t}, & w_{i,t} & \sim {\cal N}(0,\sigma_{w_i}^2), & i & =1,\ldots, n  \label{Beta:RW}
\end{align}
where $\sigma_{w_i}^2 \ge 0$, $i=1,\ldots,n$, and $\varepsilon_t$ and $w_{i,t}$ are mutually uncorrelated Gaussian white noise processes.

In summary, the most simple test for weak-form market efficiency applies equations~\eqref{AR:n:time} and \eqref{Beta:RW}, and if it yields the following condition(s): if $\beta_{1,t}=\beta_{2,t}= \ldots =\beta_{n,t}=0$ at any time $t$, then the market is weak-form efficient at that time. Under this condition, equation~\eqref{AR:n:time} simply reduces to a Gaussian white noise process for the return series, i.e. the returns are serially uncorrelated over the time period examined. Thus, the modeling goal is to estimate the hidden state-vector $[\beta_{1,t},\beta_{2,t}, \ldots, \beta_{n,t}]$. This is frequently achieved by applying the classic Kalman filter (KF) approach, as elucidated in~\cite{2009:ItoSugiyama}. We emphasize the fact that first-order autoregressive models,  AR(1), are traditionally used in practice because they are rich enough to be interesting and simple enough to permit complicated extensions for modeling the dynamics of the coefficients.

One of the most important and sophisticated extensions possible is to model the heteroscedasticity frequently observed in return series and embed it in the tests for evolving weak-form efficiency. This approach allows the modeling of a phenomenon known as volatility clustering, first observed and noted in~\cite{1967:Mandelbrot}: ``large changes tend to be followed by large changes, of either sign, and small changes tend to be followed by small changes''. Thus, in contrast to the simple test summarized by equations~\eqref{AR:n:time} and \eqref{Beta:RW} with the assumption of homoscedastic conditional variance in~\eqref{AR:n:time}, more sophisticated modeling approaches suggest using GARCH-type models in the conditional variance equation. Models that capture the time-varying nature of the volatility process and its clustering effect are referred to as Autoregressive Conditional Heteroscedasticity (ARCH) models, originally defined in~\cite{1982:Engle} and later extended to generalized ARCH (GARCH) models in~\cite{1986:Bollerslev}.

Thus, multiple tests for weak-form  market efficiency of varying degrees of complexity with implied, conditional heteroscedasticity can be conducted using an AR($n$) process~\eqref{AR:n:time} with time-varying regression coefficients $\beta_{i,t}$, $i=1,\ldots,n$ modeled by the random walk~\eqref{Beta:RW}, and with a time-varying volatility  $\sigma_t^2:=\sigma_{\varepsilon}^2(t)$ following a GARCH(p,q) process:
\begin{align}
\sigma^2_t & =  \omega + \sum \limits_{l=1}^{p} a_l \varepsilon_{t-l}^2 + \sum \limits_{j=1}^{q} b_j \sigma^2_{t-j}. \label{garch:pq}
\end{align}

Time-varying AR(1) models combined with GARCH(1,1)-type processes are the most widely used modeling techniques for measuring any degree of market efficiency. The dynamics of the {\it evolving market efficiency} is then reflected by the AR(1) coefficient $\beta_{1,t}$, $t=1, \ldots, N$, and may be recovered by the KF as demonstrated in~\cite{1997:Emerson,1999:Zalewska,1991:Hall:note}. The state estimation problem is often combined with parameter estimation by using the method of maximum likelihood, since econometric models are typically parameterized. The entire state and parameter estimation scheme is explained in detail in~\cite{2019:RJ:Kulikov}, where the classic KF should be employed over the extended KF in the case of the linear systems examined in this paper. Stable gradient-based estimation methods and their application to time-varying econometric models can be also found in~\cite{2021:KulikovaTsyganovaKulikov,2020:IEEE:KulikovaTsyganovaKulikov}.

We now present the key GARCH-type models that have been used in a variety of empirical studies that have applied the test for evolving market efficiency. First, we consider the most simple case where a GARCH(1,1) process is employed. This model was used to examine the Russian and Czech stock markets in~\cite{2002:Hall,2008:Posta}, and is given as follows:
\begin{align}
y_t & =  \beta_{1,t} y_{t-1} + \sigma_t \varepsilon_t, \quad \varepsilon_t \sim {\cal N}(0,1) \label{AR:n:time1} \\
\beta_{1,t} & = \beta_{1,t-1} + w_{1,t}, \quad  w_{1,t} \sim {\cal N}(0,\sigma_{w_1}^2) \label{Beta:RW1} \\
\sigma^2_t & = \omega + a_1 (y_{t-1} - \beta_{1,t-1} y_{t-2})^2 + b_1 \sigma^2_{t-1} = \omega + a_1 (\sigma_{t-1}\varepsilon_{t-1})^2 + b_1 \sigma^2_{t-1} \label{garch:11}
\end{align}
where $\omega>0$, $a_1 \ge 0$, $b_1 \ge 0$ and $a_1 + b_1 <1$ ensure that the GARCH(1,1) process is stationary. The disturbances $\varepsilon_t$ and $w_{i,t}$ are mutually uncorrelated (Gaussian) white noises. The unconditional expectation of the conditional variance is constant and finite, and is given by $\omega/(1-a_1-b_1)$. The test for weak-form market efficiency requires the condition $\beta_{1,t}=0$ to hold at some time $t$. This yields the reduced formula $y_t = \sigma_t\varepsilon_t$ where $\varepsilon_t \sim {\cal N}(0,1)$ is a Gaussian white noise. The returns are then serially uncorrelated.

A more sophisticated model incorporates a GARCH-in-Mean(1,1) process; see the first ARCH-in-Mean specification introduced in~\cite{1987:Engle}. Investigating time-varying market efficiency based on a GARCH-in-Mean(1,1) process was proposed in~\cite{1997:Emerson,1999:Zalewska} where it was termed the {\it test for evolving efficiency} (TEE). This approach is now the most commonly used. The TEE is given as follows:
\begin{align}
y_t & =  \beta_{1,t} y_{t-1} + \delta \sigma_t^2 + \sigma_t \varepsilon_t, \quad \varepsilon_t \sim {\cal N}(0,1) \label{TEE:yk1} \\
\beta_{1,t} & = \beta_{1,t-1} + w_{1,t}, \quad  w_{1,t} \sim {\cal N}(0,\sigma_{w_1}^2) \label{TEE:Beta} \\
\sigma^2_t & = \omega + a_1 (y_{t-1} - \delta \sigma_t^2 - \beta_{1,t-1} y_{t-2})^2 + b_1 \sigma^2_{t-1} = \omega + a_1 (\sigma_{t-1}\varepsilon_{t-1})^2 + b_1 \sigma^2_{t-1} \label{TEE:garch}
\end{align}
where the disturbances $\varepsilon_t$ and $w_{i,t}$ are mutually uncorrelated white noises and, again, the test for weak-form market efficiency requires that the condition $\beta_{1,t}=0$ holds at any time $t$. The GARCH-in-Mean(1,1) specification reduces to the classic GARCH(1,1) when $\delta = 0$.

Both models in equations~\eqref{AR:n:time1}~-- \eqref{garch:11} or in equations~\eqref{TEE:yk1}~-- \eqref{TEE:garch} account for volatility clustering, but they are not capable of modeling the asymmetric (leverage) effect observed in real return series, i.e. that bad news has a larger impact on the future volatility of asset returns than good news. This effect can be taken into account by employing asymmetric GARCH-type processes. We follow the empirical study in~\cite{2000:Rockinger} where the first test for evolving market efficiency based on an asymmetric threshold heteroscedastic model was developed as follows:
\begin{align}
y_t & =  \beta_{1,t} y_{t-1} + \sigma_t \varepsilon_t, \quad \varepsilon_t \sim {\cal N}(0,1) \label{Asym1:yk1} \\
\beta_{1,t} & = \beta_{1,t-1} + w_{1,t}, \quad  w_{1,t} \sim {\cal N}(0,\sigma_{w_1}^2) \label{Asym1:Beta} \\
\sigma^2_t = \; & \omega + a_1^+ (y_{t-1} - \beta_{1,t-1} y_{t-2})^2\Lambda_{(y_{t-1} - \beta_{1,t-1} y_{t-2})>0} \,+ \nonumber \\
& a_1^- (y_{t-1} - \beta_{1,t-1} y_{t-2})^2\Lambda_{(y_{t-1} - \beta_{1,t-1} y_{t-2})<0} +  b_1 \sigma^2_{t-1} \nonumber \\
 = \; &\omega + a_1^+ (\sigma_{t-1}\varepsilon_{t-1})^2\Lambda_{(\sigma_{t-1}\varepsilon_{t-1})>0} +a_1^- (\sigma_{t-1}\varepsilon_{t-1})^2\Lambda_{(\sigma_{t-1}\varepsilon_{t-1})<0} +  b_1 \sigma^2_{t-1} \label{Asym1:garch}
\end{align}
where the disturbances $\varepsilon_t$ and $w_{i,t}$ are mutually uncorrelated white noises and, again, the test for weak-form market efficiency  implies the condition $\beta_{1,t}=0$ at any time $t$. The threshold GARCH (T-GARCH) specification in equation~\eqref{Asym1:garch} was originally proposed in~\cite{1992:Campbell,1993:Glosten,1994:Zakoian}. Here, two dummy variables are included: $\Lambda_{(\cdot)>0}$ takes the value $1$ when the last period's error is positive, and the value $0$ otherwise. Similarly $\Lambda_{(\cdot)<0}$ takes the value $1$ if the last period's error is negative, and the value $0$ otherwise. This allows the conditional standard deviation to depend upon the sign of the lagged innovations. The threshold GARCH(1,1) specification in equation~\eqref{Asym1:garch} reduces to the classic GARCH(1,1) in equation~\eqref{garch:11} when $a_1^+ = a_1^-$, i.e. when the parameter value $a_1^+$, which measures the impact of a past positive error on the current variability of returns, coincides with $a_1^-$, which measures the impact of a past negative error on the current variability of returns.

An alternative asymmetric GARCH modeling approach was explored in~\cite{2003:Li:AMH:China}. Here, the proposed model simultaneously accounts for the leverage effect, examines the possibility of information transmission, and tests for evolving market efficiency. For this, the conditional variance of returns is assumed to evolve according to an alternative asymmetric GARCH (A-GARCH) process with an additional term $(\varepsilon^+_{t-1})^2$, and may be summarized as:
\begin{align}
y_t & =  \beta_{1,t} y_{t-1} + \sigma_t \varepsilon_t, \quad \varepsilon_t \sim {\cal N}(0,1) \label{Asym2:yk1} \\
\beta_{1,t} & = \beta_{1,t-1} + w_{1,t}, \quad  w_{1,t} \sim {\cal N}(0,\sigma_{w_1}^2) \label{Asym2:Beta} \\
\sigma^2_t & = \omega + a_1 (y_{t-1} - \beta_{1,t-1} y_{t-2})^2 + a_1^+ (u^+_{t-1})^2 + b_1 \sigma^2_{t-1} \nonumber \\
& = \omega + a_1 (\sigma_{t-1}\varepsilon_{t-1})^2 + a_1^+ (u^+_{t-1})^2 + b_1 \sigma^2_{t-1} \label{Asym2:garch}
\end{align}
where $u^+_{t-1} = \max\{(y_{t-1} - \beta_{1,t-1} y_{t-2}),0\}=\max\{\sigma_{t-1}\varepsilon_{t-1},0\}$, the disturbances $\varepsilon_t$ and $w_{i,t}$ are mutually uncorrelated white noises and, again, the test for weak-form market efficiency implies the condition $\beta_{1,t}=0$ at any time $t$. The value of the parameter $a_1^+$ only reflects the impact of a past positive error on the current variability of returns.

\subsection{Estimation results for evolving efficiency of BRICS markets and discussion}

In our empirical study, we estimate all five models above using the method of maximum likelihood and summarize the resulting Akaike information criterion (AIC) values in Table~\ref{tab:aics}. The model with the lowest AIC value is preferable. To ensure a robust comparison, all models are estimated using the same filtering method (the classic KF), the same parameter estimation strategy (the method of maximum likelihood), the same optimization method (the built-in optimizer {\tt fmincon} in MATLAB), the same initial parameter values, and the same initialization of the filter. More details about estimation procedure and the state-space methodology utilized can be found in Appendix.

\begin{table}[h]
{\scriptsize
\caption{The AIC values of the estimated time-varying AR(1) models with homoscedastic and heteroscedastic assumptions. The lowest AIC denotes the preferred model by $\blacktriangledown$ in each column.} \label{tab:aics}                                                                                                                                                                                                                                                                                                                                                                                                                                                                                                                                                                                                                                                                \begin{tabular}{llllllll}
\toprule
{\bf Model} & {\bf Effect(s) modeled} & {\bf Ref.} & {\bf Brazil} & {\bf Russia} &  {\bf India} & {\bf China} & {\bf S. Africa}\\
\toprule
Time-varying AR(1) & Homoscedasticity       & in~\cite{2009:ItoSugiyama}          & -1205.71 & -929.81  & -1343.86  & -1262.33 & -1449.92 \\
GARCH(1,1)         & Heteroscedasticity      & in~\cite{2002:Hall}                & -1225.34 & -1037.95$\blacktriangledown$
                                                                                                         & -1360.70$\blacktriangledown$
                                                                                                                     & -1299.63 & -1496.38 \\
T-GARCH(1,1)   & Heteroscedasticity, leverage & in~\cite{2000:Rockinger}      & -1223.85 & -1036.36 & -1360.11  & -1302.74$\blacktriangledown$ & -1506.19$\blacktriangledown$ \\
A-GARCH(1,1)  & Heteroscedasticity, leverage & in~\cite{2003:Li:AMH:China}      & -1233.25 & -1035.95 & -1358.70  & -1301.71 & -1494.38 \\
GARCH-M(1,1) &Heteroscedasticity  & in~\cite{1999:Zalewska} & -1234.21$\blacktriangledown$
                                                                                              & -1037.13 &  -1358.86 & -1297.69 & -1492.29 \\
\bottomrule
\end{tabular}}
\end{table}

As can be seen from the first row of Table~\ref{tab:aics}, the tests for time-varying market efficiency assuming a homoscedastic variance are rejected for all five indices. This model is too restrictive when compared to the models that assume heteroscedasticity and GARCH-type processes. Next, it is interesting to note that the asymmetric T-GARCH specification in~\eqref{Asym1:garch} provides the lowest AIC values for the Chinese and South African series. This is a strong indication of the existence of a leverage effect for these stock markets. Significant asymmetry here may be due to government intervention in investors' activity. The empirical study in~\cite{2003:Li:AMH:China} revealed clear evidence of a leverage effect in the Shanghai index, which corresponds well with our results here. Empirical studies of the South African market dynamics have never reported any evidence of asymmetry. However, it should be noted that the most recent empirical study~\cite{2018:Heymans} only considered data ending in March 2015. This period excludes a critical period of political instability that began in 2015 and resulted in dramatic changes in the level of market efficiency in South Africa. Finally, it is clear that for Russian and Indian stock markets the symmetric GARCH specification provides the best fit for evolving market efficiency. This is a sign of an insignificant leverage effect in Russian and Indian markets. For the Brazilian market, the best-fit model is the TEE test based on the GARCH-in-Mean specification, which also implies an insignificant leverage effect.

\begin{table}[h]
{\scriptsize
\caption{Estimation results of the best-fit, time-varying AR(1) models with GARCH-type processes for the time-varying conditional variance. Standard errors are given in parentheses.} \label{tab:est}                                                                                                                                                                                                                                                                                                                                                                                                                                                                                                                                                                                                                                                                \begin{tabular}{lcccccc}
\toprule
 &  {\bf Brazil} & {\bf Russia} &  {\bf India} & {\bf China} & {\bf S. Africa}\\
\toprule
Best-fit specification: & {\tt GARCH-in-Mean(1,1)} & {\tt GARCH(1,1)} & {\tt GARCH(1,1)}  & {\tt T-GARCH(1,1)} & {\tt T-GARCH(1,1)}  \\
Model summarized by: & eqs.~\eqref{TEE:yk1}~-- \eqref{TEE:garch} & eqs.~\eqref{AR:n:time1}~-- \eqref{garch:11} & eqs.~\eqref{AR:n:time1}~-- \eqref{garch:11}  & eqs.~\eqref{Asym1:yk1}~-- \eqref{Asym1:garch} & eqs.~\eqref{Asym1:yk1}~-- \eqref{Asym1:garch} \\
\hline
\\
$\hat \sigma_{w_1}^2$ & 0.00001   & 0.00001   & 0.00001    & 0.00001   &  0.00002 \\
                      & (0.00022) & (0.00000) & (0.00063)  & (0.00000) &  (0.00016) \\[3pt]
$\hat \omega$         & 0.00078   & 0.00055   & 0.00022    & 0.00037 &  0.00037 \\
                      & (0.00038) & (0.00025) & (0.00014)  & (0.00015) &  (0.00018) \\[3pt]
$\hat \delta$         & -0.06173  & ---       & ---        & --- &  --- \\
                      & (0.09301) &           &            &     &      \\[3pt]
$\hat a_1$            & 0.27552   & 0.25102	  & 0.08657    & --- &  --- \\
                      & (0.11273) & (0.05356) & (0.04092)  &  &      \\[3pt]
$\hat b_1$            & 0.65540   & 0.73666   & 0.86794    & 0.78901 &  0.67738 \\
                      & (0.10894) & (0.05488) & (0.06067)  & (0.06075) &  (0.10519) \\[3pt]
$\hat a_1^+$          & ---       & ---       &  ---       & 0.19650 &  0.04515 \\
                      &           &           &            & (0.05931) &  (0.06544) \\[3pt]
$\hat a_1^-$          & ---       & ---       & ---        & 0.09102 &  0.32262 \\
                      &           &           &            & (0.05071) &  (0.09056) \\
\hline
max log LF            & 622.1077  & 522.9770 & 684.3517  & 656.3726 &  758.0979 \\
AIC                   & -1234.21  & -1037.95 & -1360.70  & -1302.74 & -1506.19 \\
\bottomrule
\end{tabular}}
\end{table}

In Table~\ref{tab:est}, we summarize the parameter estimates obtained by the method of maximum likelihood for each best-fit model and return series under examination. The most interesting cases are the Chinese and South African stock markets, because of the aforementioned leverage effect. Having examined the estimates for the South African market, we conclude that a strong asymmetric effect is present because the large difference between the estimates for $\hat a_1^+$ and $\hat a_1^-$. In addition, the impact of past negative errors on the current variability of returns is much higher than the impact of past positive errors, i.e. the market is mainly affected by negative news. In contrast, the Chinese market exhibits a weaker leverage effect and the estimates for $\hat a_1^+$ and $\hat a_1^-$ are quite close to each other. However, the impact of past positive errors is slightly higher than the negative errors.
The GARCH-in-Mean specification in the TEE methodology provides the best-fit model for the Brazilian market. The estimated value of the parameter $\hat \delta$ is not far from zero and, hence, the estimated model is close to the classic GARCH(1,1) specification.
Finally, the symmetric GARCH(1,1) specification provides the best-fit model for the Russian and Indian markets. The estimates substantiate the stability of the models since $\hat a_1 + \hat b_1 <1$ and $\hat a_1 << \hat b_1$ in each.

\begin{figure}[th!]
\begin{tabular}{cc}
\includegraphics[width = 0.45\textwidth]{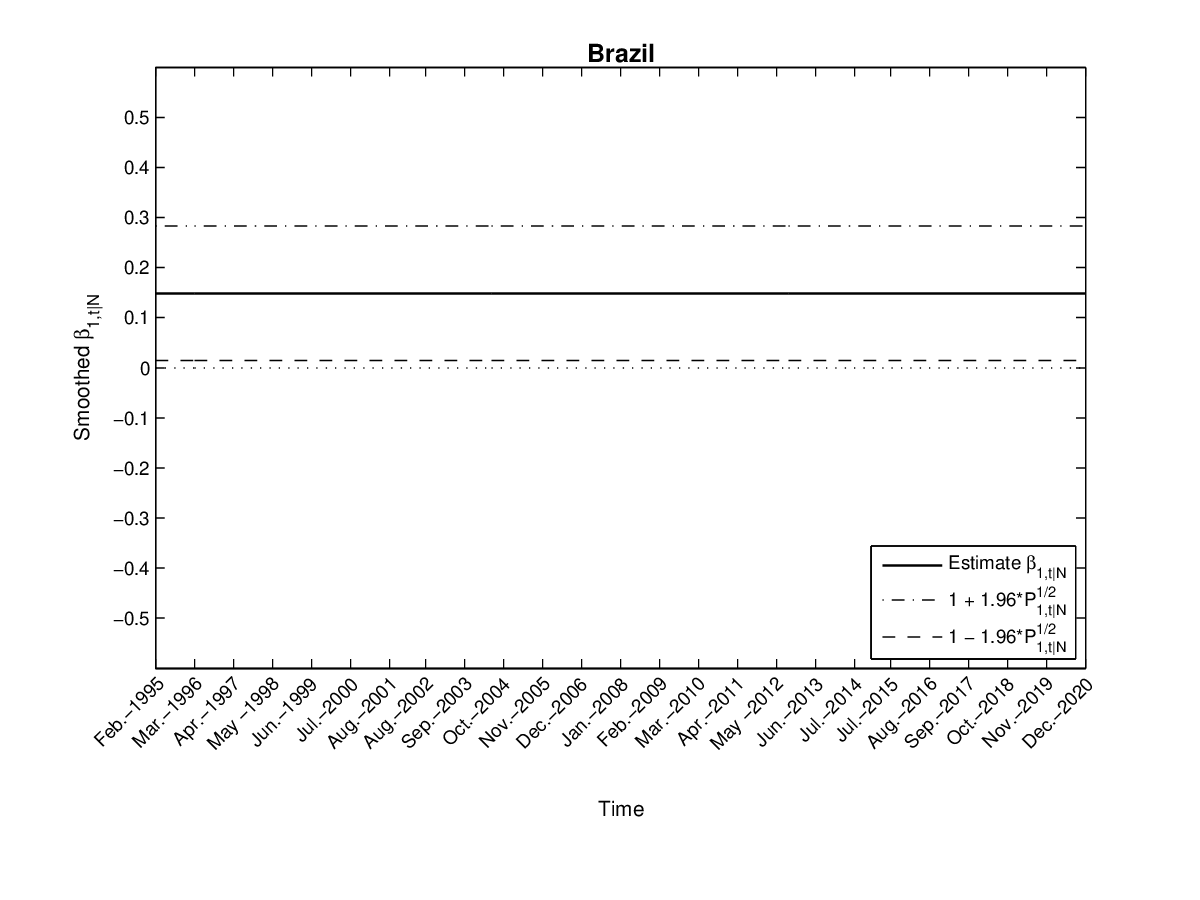} &
\includegraphics[width = 0.45\textwidth]{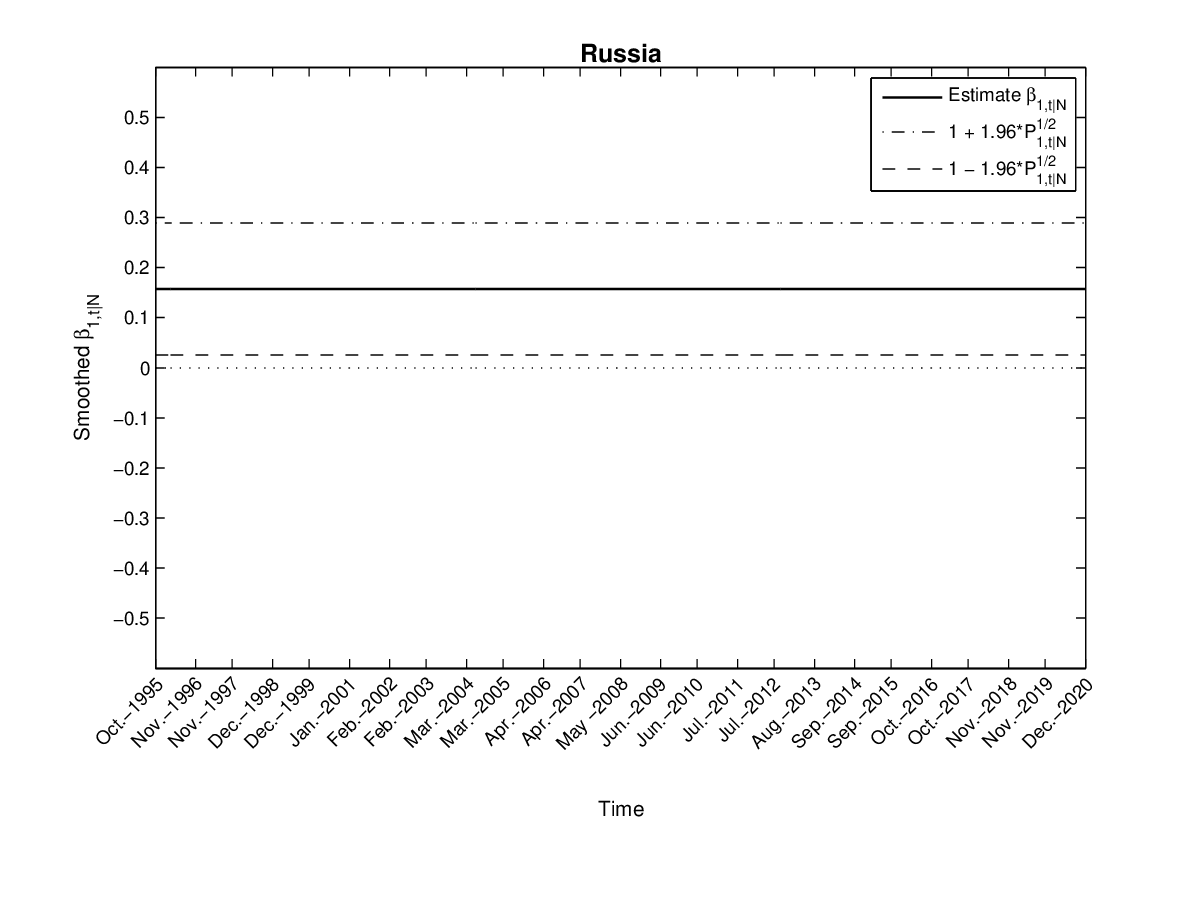} \\
\includegraphics[width = 0.45\textwidth]{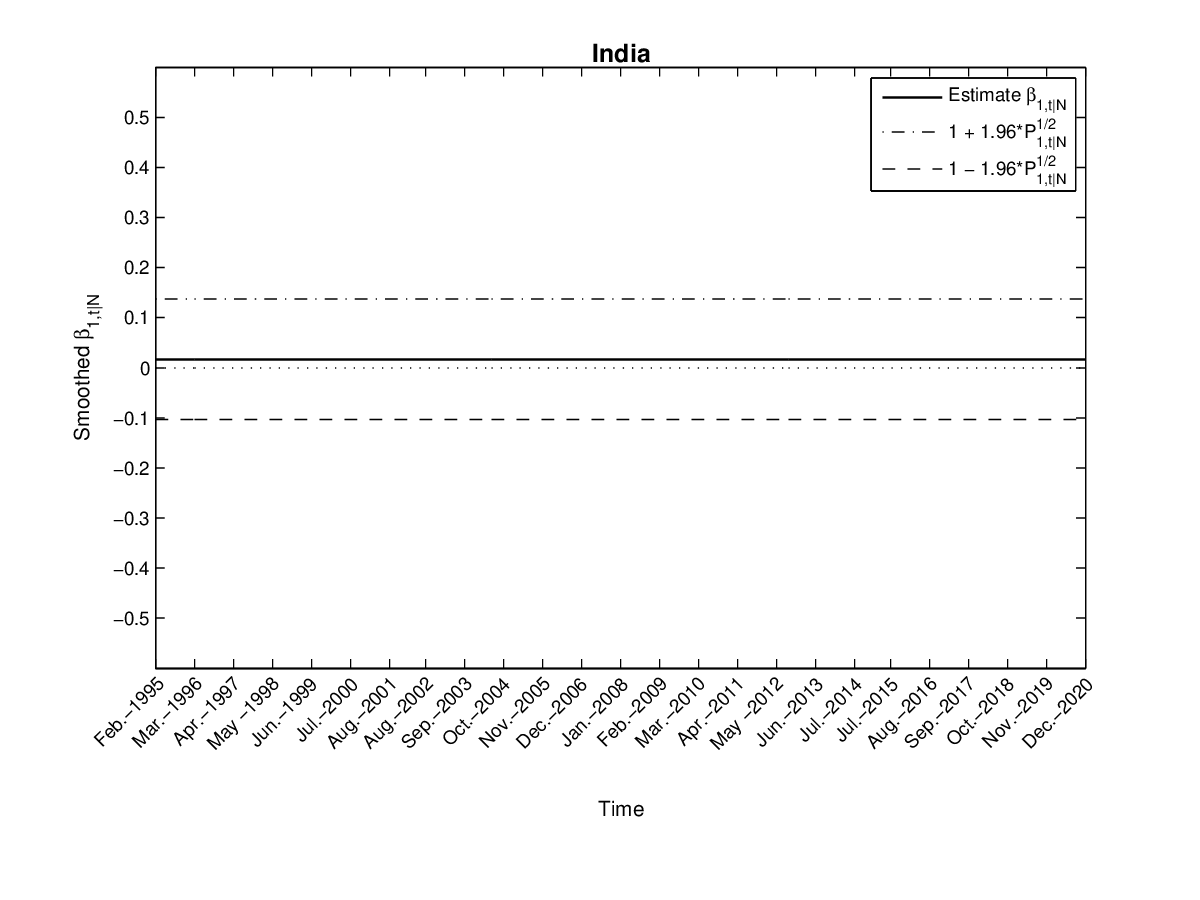} &
\includegraphics[width = 0.45\textwidth]{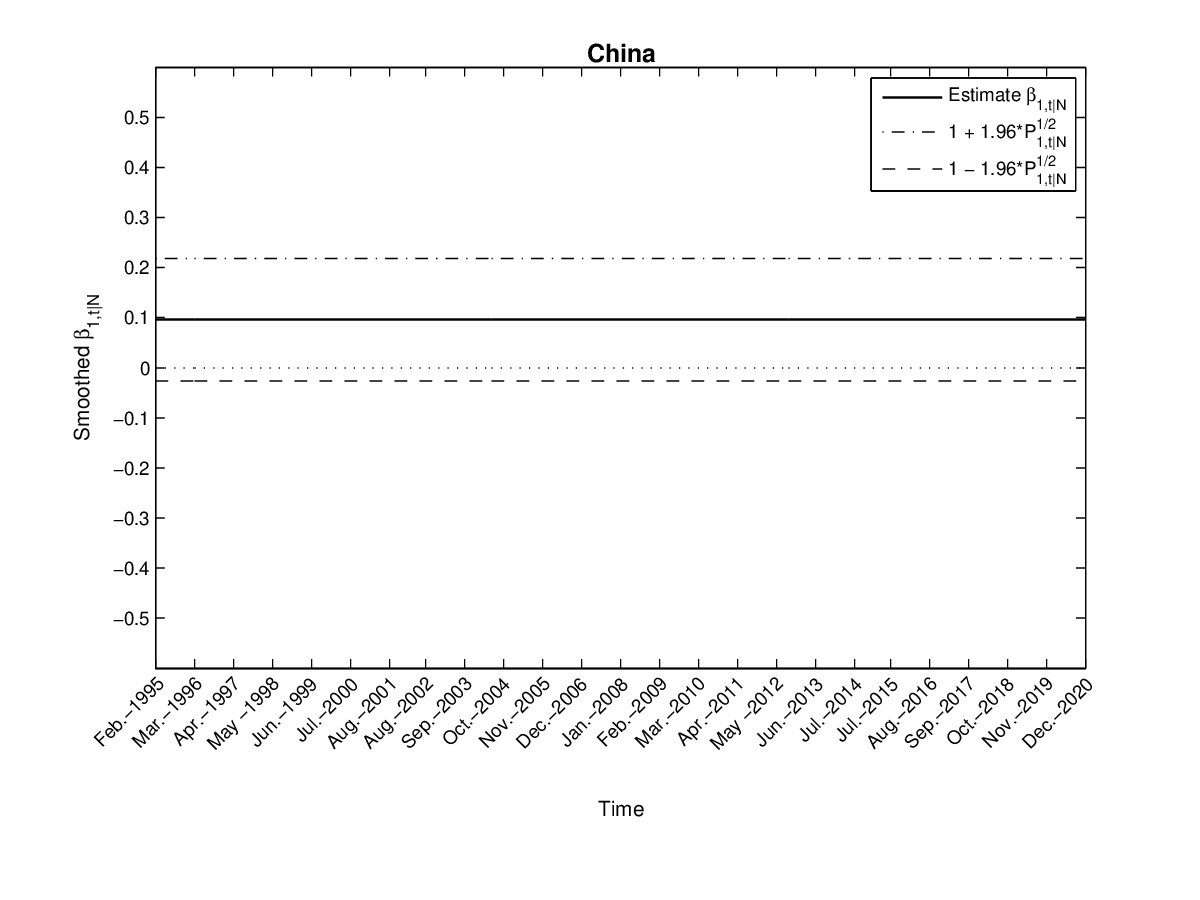} \\
\includegraphics[width = 0.45\textwidth]{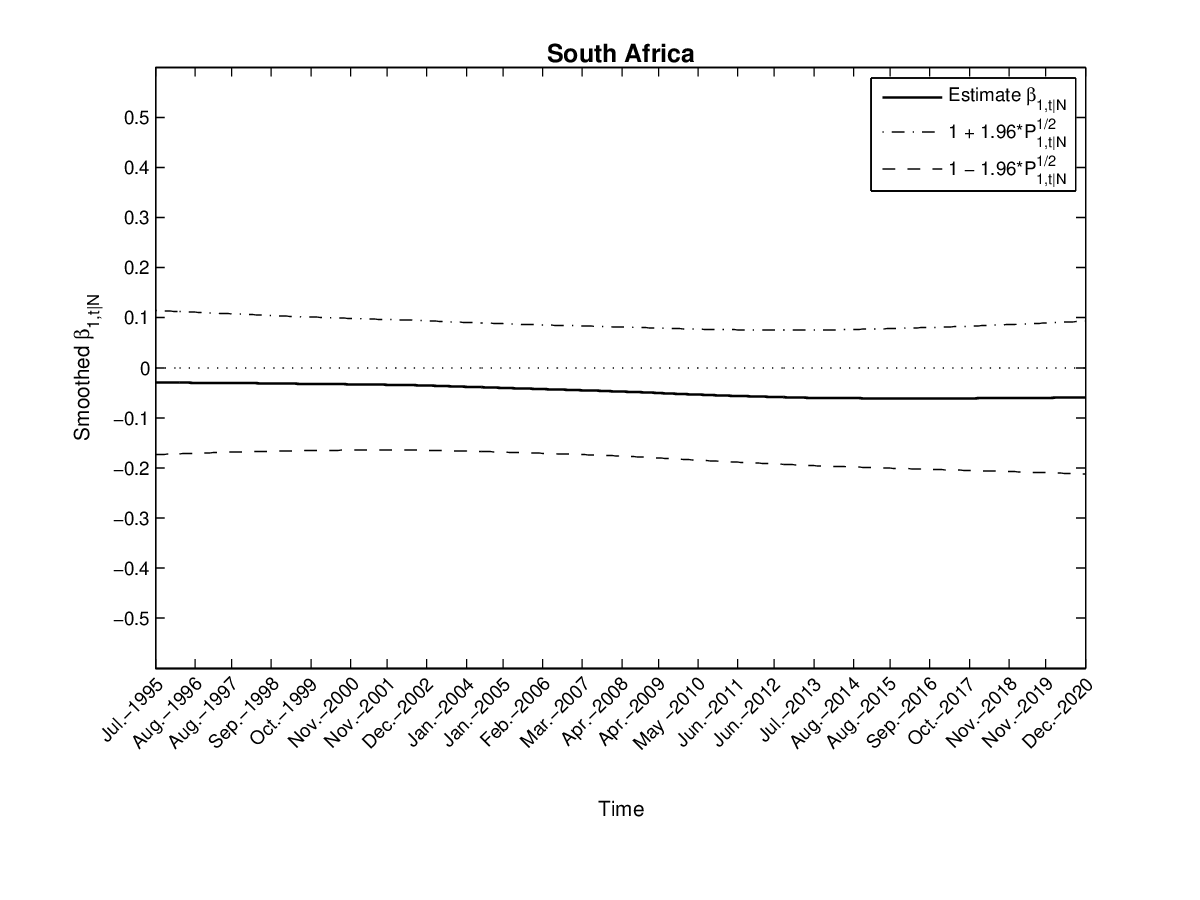} &
\includegraphics[width = 0.45\textwidth]{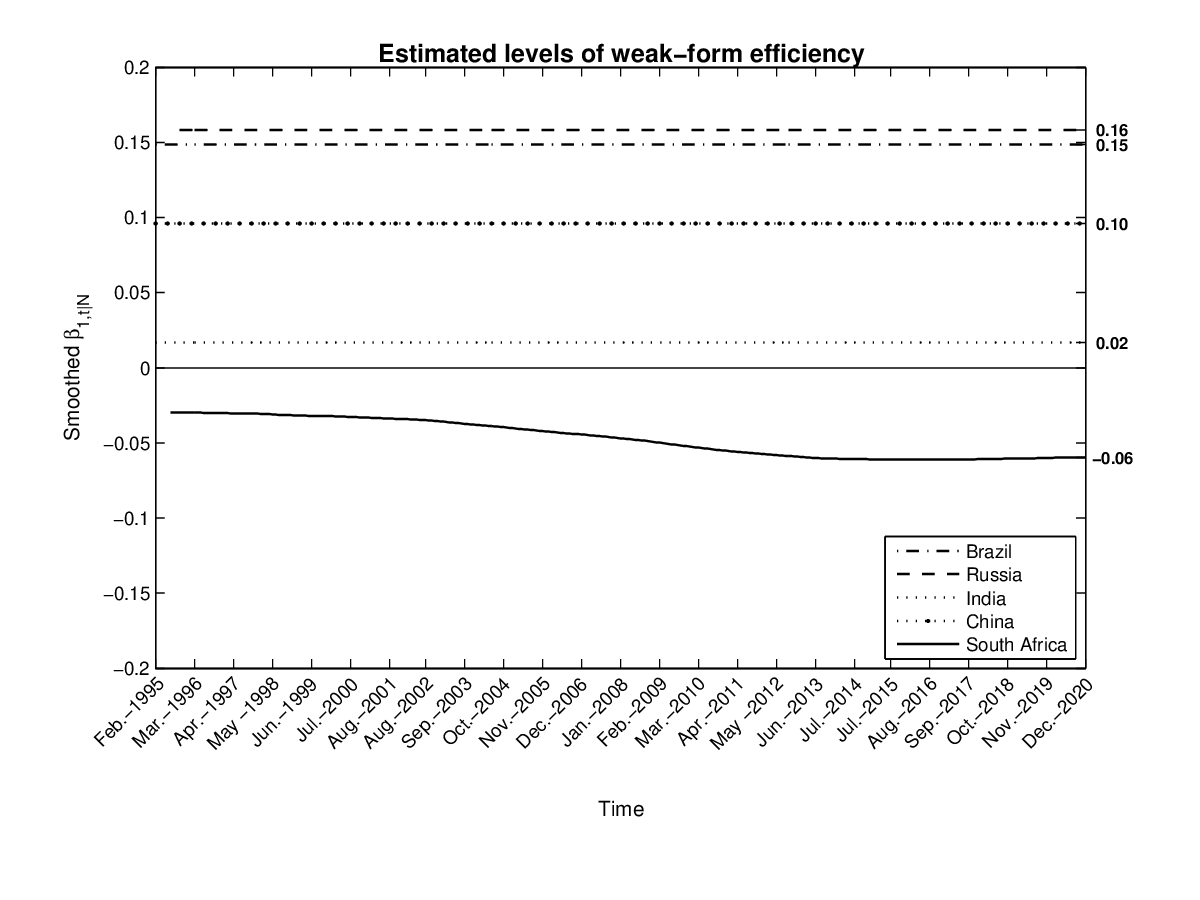}
\end{tabular}
\caption{Smoothed estimates of the weak-form market efficiency levels $\hat \beta_{1,t}$ recovered from the best-fit models for each market under examination; see Table~\ref{tab:aics}. The optimal parameters values of the best-fit model calibration are summarized in Table~\ref{tab:est}. The bottom right graph illustrates the time-varying $\beta_{1,t}$ estimations for all series for comparative reasons.}
\label{fig:TEE}
\end{figure}

In Figure~\ref{fig:TEE}, we provide plots of the time-varying $\beta_{1,t}$ coefficients for each best-fit model and series. Having compared the results, we conclude the following. As can be seen in the last graph, the Russian stock market is the least efficient, with an estimated  value of $\hat \beta_{1,t} = 0.16$ at the end of our study on December 2020. The estimated degrees of (in)efficiency for Brazil and Russia are very similar. The most stable evolving efficiency process with the smallest degree of inefficiency is India, which has an estimated value of $\hat \beta_{1,t} = 0.02$ at the end of our study. The South African market was weak-form efficient at the beginning of our empirical study, with an estimated level of efficiency very close to that of India's. However, it gradually becomes less efficient over the course of the study. Despite this, it remains the second most efficient country among the BRICS group.

Finally, we note that our modeling technique, which applies the classic KF/smoother to calibrate the models and estimate the hidden efficiency process, only recovers the main trend $\hat \beta_{1,t}$ in the changing degree of a weak-form efficiency.

\section{Concluding remarks} \label{Section:conclusion}

In this paper adaptive, weak-form efficiency changes in the BRICS stock markets have been estimated over a period that includes the 2008-2009 global financial crisis and the recent COVID-19 recession. We applied rolling window tests for sample autocorrelation as well as appropriate time-varying autoregressive models with both homoscedastic and heteroscedastic conditional variance assumptions. The models were calibrated by the method of maximum likelihood and used the Kalman filter to recover the hidden dynamics of the market efficiency process. The rolling window tests reveal that all BRICS stock markets exhibit a trend toward weak-form efficiency, except China, which displays the opposite trend and gradually becomes less efficient in the last few years. China has also been severley affected by the recent COVID-19 recession. The Indian and South African markets experienced significant turbulence in 2019-2020, but they show a strong recovery process in the last year. Brazilian and Russian markets were less impacted by the COVID-19 recession, but they were greatly affected by the 2008-2009 global financial crisis. Finally, GARCH-type modeling techniques have substantiated some conclusions about the different degrees of time-varying efficiency of the BRICS markets. We found that the Russian stock market is the least efficient and its estimated degree of (in)efficiency is very close to the Brazilian case. The most stable evolving efficiency process is displayed by India with the smallest degree of inefficiency estimated at the end of our empirical study. We have also detected a strong sign of the leverage effect in the South African market that might be explained by its government's intervention in investors' activity. China exhibits a much weaker leverage effect, but an asymmetric impact of past errors on the current variability of returns has been detected.

The state-space approach alongside a classical KF that we have applied here is conventionally used in econometric literature. This methodology is limited by the linearity of the classical KF, which may reduce its estimation and forecasting qualities since the proposed autoregressive models (with time-varying coefficients and heteroscedastic variance) can be regarded as nonlinear stochastic systems affected by multiplicative state noises. Thus, our future plans concern the derivation of nonlinear Bayesian filters for more accurate tracking of the dynamics of weak-form market efficiency. Additional interesting topics for future research include the following: i) a more delicate modeling of the dynamics of the regression coefficients and the changes in conditional variance over time, and ii) improving the quality of tracking by designing higher-order nonlinear estimators where special attention can be paid to derivative-free filtering techniques, which simplifies their practical application.

\section*{Acknowledgements}
The first and third authors acknowledge the financial support of the Portuguese FCT~--- \emph{Funda\c{c}\~ao para a Ci\^encia e a Tecnologia}, through the projects UIDB/04621/2020 and UIDP/04621/2020 of CEMAT/IST-ID, Center for Computational and Stochastic Mathematics, Instituto Superior T\'ecnico, University of Lisbon and through the \emph{Scientific Employment Stimulus - 4th Edition} (CEEC-IND-4th edition) programme, grant number 2021.01450.CEECIND.

\section*{References}
\bibliographystyle{model1b-num-names}
\bibliography{BibTex_Library/books,%
              BibTex_Library/AME_Tests,%
              BibTex_Library/Lit_Finance}

\appendix
\section{Classical Kalman filter and state-space model calibration} \label{sec:appendix}

The estimation problem considered in this paper is that of the unknown dynamic state and system parameters of linear discrete-time system written in the state-space form as follows:
 \begin{align}
   x_{t} = & F(\theta) x_{t-1} + B(\theta) u_{t-1} + G(\theta) w_{t},  & w_t & \sim {\cal N}\left(0, Q(\theta)\right), \quad Q(\theta) \ge 0 \label{eq:state:space1}
 \end{align}
where $x_t:= x(t_k)$ is an unknown $n$-vector to be estimated, $u_t$ is a control $d$-vector, and $\{w_t, t=1,\ldots \}$ is a white Gaussian sequence $q$-vector. The matrices are $F \in {\mathbb R}^{n\times n}$, $B \in {\mathbb R}^{n\times d}$, $G \in {\mathbb R}^{n\times q}$ and $Q \in {\mathbb R}^{q\times q}$. The distribution of the initial condition $x_0$ is assumed given, say $x_0 \sim {\cal N}(\bar x_0, P_0)$, and $x_0$ is independent of $\{w_t\}$. Additionally, the state-space model is parameterized by the unknown system parameter vector $\theta \in {\mathbb R}^{p}$, which should be estimated together with the hidden state vector $x_t$ from the noisy data at hand.

The noisy $m$-vector observations (measurements) $y_t:=y(t_k)$ are given by
 \begin{align}
   y_{t} = & H(\theta) x_t + \varepsilon_t, & \varepsilon_t & \sim {\cal N}\left(0,R(\theta)\right), \quad R(\theta) >0 \label{eq:state:space2}
 \end{align}
where $H \in {\mathbb R}^{m\times n}$ and $\{\varepsilon_t: t=1,\ldots \}$ is an $m$-vector, Gaussian white sequence where $R \in {\mathbb R}^{m\times m}$. The processes $\{ w_t\}$ and $\{\varepsilon_t \}$ are assumed to be independent, and $\{ \varepsilon_t \}$ is independent of the initial $x_0$.

The goal of any filtering technique is to recover an unknown (hidden) random process $X = \{x_t\}_{t=0}^{T} = \{x_0, \ldots, x_T\}$ from an observed
sequence $Y = \{y_t\}_{t=0}^{T} = \{y_0, \ldots, y_T\}$.  Since model~\eqref{eq:state:space1}, \eqref{eq:state:space2} is also characterized by an additional, unknown set of system parameters $\theta$, both the dynamic state and the unknown system parameters should be estimated from only the observed (noisy) signal. This problem can be resolved by {\it adaptive} filtering techniques where both the states $X = \{x_t\}_{k=0}^{T} = \{x_0, \ldots, x_T\}$ and the system parameters $\theta$  are simultaneously estimated. Initially, we should choose a performance index (cost function) that reflects the difference between the actual system and the filter associated with the model. The method of maximum likelihood is a general method that can be successfully applied here. The resulting maximum likelihood estimate $\hat \theta^*$ has some attractive properties. It is an asymptotically unbiased, consistent and efficient estimate. For state-space models, the log likelihood function (LF) can be evaluated, as explained in~\cite{1989:Harvey:book}:
\begin{align}
\ln L\left(\theta | Y_0^T\right)  = - c_0 - \frac{1}{2}\sum \limits_{t=0}^T \left\{\ln\left(\det R_{e,t}\right) + e_t^\top R_{e,t}^{-1}e_t \right\} \label{eq:llf}
\end{align}
with $c_0$ being a constant and the residual $e_t$ and its covariance matrix $R_{e,t}$, $t=0,\ldots, T$, are calculated from the Kalman filtering (KF) equations, e.g. see~\cite[Chapter~9]{2000:Kailath:book}. This method for calculating the log LF is often referred to
as prediction error decomposition in econometric literature, see~\cite{1989:Harvey:book}.

\begin{table}[ht!]
{\small
\caption{A general scheme of adaptive filtering using the classical KF.}
\begin{tabular}{p{0.45\textwidth}|p{0.5\textwidth}}
\hline
{\bf Algorithm~1} \textsc{(Adaptive scheme)} & {\bf Algorithm~2} \textsc{(classical KF)} \\
\underline{Input Data:} Initial value $\theta_{0}$. & \underline{Input Data:} Initial $\hat x_{0|-1}$, $P_{0|-1}$ and current $\theta_{n}$. \\
\underline{Process:} Substitute current approximation $\theta_{n}$ into the state-space model. Process available measurements $\{ y_0, \ldots, y_T\}$ as follows:
\begin{itemize}
\item Set the initial values $\hat x_{0|-1}=\bar x_0$, $P_{0|-1} = P_0$ of the filter to be used.
\item Apply the underlying filtering method (e.g. KF given by Algorithm~2) for current approximation $\theta_{n}$. This allows calculation of the cost function at the point $\theta_{n}$, e.g. the log LF in~\eqref{eq:llf}, and for recovering the state estimates $\{\hat x_{t|t}\}_0^T$.
\item Use the chosen optimization method to find the next point $\theta_{n+1}$.
\end{itemize}
Repeat the process for the approximation $\theta_{n+1}$ ($n=0, 1, \ldots$) until the chosen stopping criterion is satisfied.
&
\underline{Process:} Kalman filtering

\vspace{0.2cm}
\texttt{Measurement Updates ($t=0, \ldots, T$)}. Compute the estimate $\hat x_{t|t}$ and error covariance matrix $P_{t|t} = \E{(x_t-\hat x_{t|t})(x_t-\hat x_{t|t})^{\top}}$  as follows:
\begin{eqnarray}
e_t          & =  & y_t-H_t\hat x_{t|t-1}, \\
R_{e,t}      & = & H_tP_{t|t-1}H_t^{\top} + R_t \\
K_t          & = & P_{t|t-1}H_t^{\top} R_{e,t}^{-1} \\
\hat x_{t|t} & = & \hat x_{t|t-1}+K_t e_t \\
P_{t|t}      & = & (I -K_tH_t)P_{t|t-1}
\end{eqnarray}
\texttt{Time Updates}. Propagate the estimate and covariance as follows:
\begin{eqnarray}
\hat x_{t+1|t} &  = & F\hat x_{t|t} + B u_{t} \\
P_{t+1|t}   & =  & FP_{t|t}F^{\top} + GQG^{\top}
\end{eqnarray} \\[-15pt]
\underline{Output Data:} Optimal system parameters $\hat \theta^*$ and state estimates $\{ \hat x_{t|t}\}_0^T$ from the filter that is tuned to the resulting optimal $\hat \theta^*$. One-step-ahead predicted estimates $\{ \hat x_{t+1|t}\}_0^T$ are also available from the filter.
& \underline{Output Data:} Log LF in equation~\eqref{eq:llf}. \\
\hline
\end{tabular}}
\label{table-adaptive}
\end{table}

In summary, any adaptive estimator consists of the two parts summarized in Table~1: (i) the filtering method for estimating the unknown state vector $\hat x_{t|t}$, and (ii) the optimization method for finding the $\hat \theta^*$ that maximizes the performance index, e.g. the log LF. In practice, these two parts are implemented in parallel: running the classical KF at each iteration step with respect to the unknown system parameters $\theta$ to generate $\{ e_t, R_{e,t} \}$ for the log LF evaluation~\eqref{eq:llf} and optimizing the log LF. In our empirical study we use the built-in MATLAB function \verb"fmincon" for optimization purposes with the default stopping criterion $\|\theta_{n} -\theta_{n-1} \| \le 10^{-6}$. The state initial values to be utilized by the KF are set to $\bar x = 0$ and $P_0 = I$ where $I$ is an identity matrix of a proper size.

Following~\cite{1989:Harvey:book}, if a model can be put in state space form, a number of powerful statistical results immediately become available such as optimal prediction via the KF, optimal estimates of unobserved components via smoothing, and estimation of unknown model parameters in the
models. In the rest of the paper, we explicitly write down the state-space form for all models used in our research.

\section{Linear state-space form for time-varying AR($n$) models}

In econometric discipline, the time-varying AR($n$) models presented above are traditionally estimated by the classical KF. For that, system~\eqref{AR:n:time}, \eqref{Beta:RW} is casted in the state-space form where the measurement matrix is replaced by a history of returns as follows  (see~\cite{1997:Emerson,1991:Hall:note,2009:ItoSugiyama}, and others):
\begin{align*}
\underbrace{
\begin{bmatrix}
\beta_{1,t}\\
\beta_{2,t} \\
\vdots \\
\beta_{n,t}
\end{bmatrix}
}_{x_t}
& =
\underbrace{\begin{bmatrix}
1 & 0 & \ldots & 0  \\0 & 1 & \ldots & 0  \\
\vdots & \vdots & \ddots & \vdots  \\
0 & 0 & \ldots & 1
\end{bmatrix}
}_{F}
\underbrace{
\begin{bmatrix}
\beta_{1,t-1}\\
\beta_{2,t-1} \\
\vdots \\
\beta_{n,t-1}
\end{bmatrix}
}_{x_{t-1}}
+
\underbrace{
\begin{bmatrix}
w_{1,t}\\
w_{2,t} \\
\vdots \\
w_{n,t}
\end{bmatrix}
}_{w_{t}}, & w_t  & \sim {\mathcal N}(0,Q) & \mbox{where } & Q = \begin{bmatrix}
\sigma_{w_1}^2 & 0 & \ldots & 0  \\
0 & \sigma_{w_2}^2 & \ldots & 0  \\
\vdots & \vdots & \ddots & \vdots  \\
0 & 0 & \ldots & \sigma_{w_n}^2
\end{bmatrix} \\
y_t & =
\underbrace{
\begin{bmatrix}
y_{t-1} & y_{t-2} & \ldots & y_{t-n}
\end{bmatrix}
}_{H_t}
\begin{bmatrix}
\beta_{1,t}\\
\beta_{2,t} \\
\vdots \\
\beta_{n,t}
\end{bmatrix}
+ \varepsilon_t, & \varepsilon_t & \sim {\mathcal N}(0,R)  & \mbox{where } & R = [\sigma_{\varepsilon}^2]
\end{align*}
where $[\beta_{1,t},\beta_{2,t}, \ldots, \beta_{n,t}]^{\top}$ is the state vector to be estimated from the data, $Y_0^T:=\{y_t\}_1^T$. The unknown system parameter vector is $\theta = [\sigma_{w_1}^2, \ldots, \sigma_{w_n}^2,\sigma_{\varepsilon}^2]$.

\section{Linear state-space form for GARCH-type model}

Following~\cite{1997:Emerson,1991:Hall:note,2021:KulikovaTsyganovaKulikov,2020:IEEE:KulikovaTsyganovaKulikov}, the linear state-space form for GARCH-type models used in our research can be written as follows. Let us start with GARCH(1,1) and equations~\eqref{AR:n:time1}~-- \eqref{garch:11}. Having denoted $h_t := \sigma_t^2$ and the state vector $x_t := [h_t,\beta_{1,t}]^{\top}$, we apply the approach suggested in the cited papers and obtain the state-space form as follows:
\begin{align*}
\underbrace{
\begin{bmatrix}
h_t \\
\beta_{1,t}
\end{bmatrix}
}_{x_t}  & =  \underbrace{
\begin{bmatrix}
b_1 & 0 \\
0 & 1
\end{bmatrix}}_{F}
\underbrace{
\begin{bmatrix}
h_{t-1} \\
\beta_{1,t-1}
\end{bmatrix}
}_{x_{t-1}}
 +
\underbrace{
	\begin{bmatrix}
		\omega  & a_1 \\
        0 & 0
	\end{bmatrix}}_{B}
\underbrace{
	\begin{bmatrix}
		1 \\
       \hat e_{t-1}^2
	\end{bmatrix}}_{u_{t-1}} + 	\underbrace{\begin{bmatrix}
        0\\
        1
    \end{bmatrix}}_{G}
	\begin{bmatrix}
        w_{t}
    \end{bmatrix}, & w_t  & \sim {\mathcal N}(0,Q) & \mbox{where } & Q = [\sigma_{w_1}^2], \\
    y_t  &  = 	\underbrace{\begin{bmatrix}
        0 &  y_{t-1}
    \end{bmatrix}}_{H_t}\begin{bmatrix}
h_t \\
\beta_{1,t}
\end{bmatrix}  + \varepsilon_t, & \varepsilon_t & \sim {\mathcal N}(0,R) & \mbox{where } & R = [\hat h_{t|t-1}]
\end{align*}
where the control input is defined by $u_{t-1} := [1, \hat e_{t-1}^2]^\top$ and $\hat e_{t}$ is the KF residual, i.e. $\hat e_{t} := y_t - \hat \beta_{1,t|t-1} y_{t-1}$ at time $t$; see more detail in~\cite{1991:Hall:note}. The unknown state vector $[h_t,\beta_{1,t}]^{\top}$ and system parameters $\theta = [\omega,a_1,b_1, \sigma_{w_1}^2]$ are to be estimated from the data, $Y_0^T:=\{y_t\}_1^T$ by using the methodology presented in~\ref{sec:appendix}.

The GARCH-in-Mean(1,1) model summarized by equations~\eqref{TEE:yk1}~-- \eqref{TEE:garch} can be represented in the state-space form in a similar way, i.e. we obtain
\begin{align*}
\underbrace{
\begin{bmatrix}
h_t \\
\beta_{1,t}
\end{bmatrix}
}_{x_t}  & =  \underbrace{
\begin{bmatrix}
b_1 & 0 \\
0 & 1
\end{bmatrix}}_{F}
\underbrace{
\begin{bmatrix}
h_{t-1} \\
\beta_{1,t-1}
\end{bmatrix}
}_{x_{t-1}}
 +
\underbrace{
	\begin{bmatrix}
		\omega  & a_1 \\
        0 & 0
	\end{bmatrix}}_{B}
\underbrace{
	\begin{bmatrix}
		1 \\
       \hat e_{t-1}^2
	\end{bmatrix}}_{u_{t-1}} + 	\underbrace{\begin{bmatrix}
        0\\
        1
    \end{bmatrix}}_{G}
	\begin{bmatrix}
        w_{t}
    \end{bmatrix}, & w_t  & \sim {\mathcal N}(0,Q) & \mbox{where } & Q = [\sigma_{w_1}^2], \\
    y_t  &  = 	\underbrace{\begin{bmatrix}
        \delta &  y_{t-1}
    \end{bmatrix}}_{H_t}\begin{bmatrix}
h_t \\
\beta_{1,t}
\end{bmatrix}  + \varepsilon_t, & \varepsilon_t & \sim {\mathcal N}(0,R) & \mbox{where } & R = [\hat h_{t|t-1}]
\end{align*}
where the control input is defined by $u_{t-1} := [1, \hat e_{t-1}^2]^\top$ and $\hat e_{t}$ is the KF residual, i.e. $\hat e_{t} := y_t - \hat \beta_{1,t|t-1} y_{t-1} - \delta \hat h_{t|t-1}$ at time $t$; see more detail in~\cite{1991:Hall:note}. The unknown state vector $[h_t,\beta_{1,t}]^{\top}$ and system parameters $\theta = [\omega,a_1,b_1, \delta, \sigma_{w_1}^2]$ are to be estimated from the data, $Y_0^T:=\{y_t\}_1^T$ by using the methodology presented in~\ref{sec:appendix}.

The T-GARCH(1,1) specification in equations~\eqref{Asym1:yk1}~-- \eqref{Asym1:garch} has the following state-space form:
\begin{align*}
\underbrace{
\begin{bmatrix}
h_t \\
\beta_{1,t}
\end{bmatrix}
}_{x_t}  & =  \underbrace{
\begin{bmatrix}
b_1 & 0 \\
0 & 1
\end{bmatrix}}_{F}
\underbrace{
\begin{bmatrix}
h_{t-1} \\
\beta_{1,t-1}
\end{bmatrix}
}_{x_{t-1}}
 +
\underbrace{
	\begin{bmatrix}
		\omega  & a_1^{+} & a_1^{-}\\
        0 & 0 & 0
	\end{bmatrix}}_{B}
\underbrace{
	\begin{bmatrix}
		1 \\
       \hat e_{t-1}^2\Lambda_{(\hat e_{t-1})>0}\\
       \hat e_{t-1}^2\Lambda_{(\hat e_{t-1})<0}
	\end{bmatrix}}_{u_{t-1}} + 	\underbrace{\begin{bmatrix}
        0\\
        1
    \end{bmatrix}}_{G}
	\begin{bmatrix}
        w_{t}
    \end{bmatrix}, & w_t  & \sim {\mathcal N}(0,Q), Q = [\sigma_{w_1}^2], \\
    y_t  &  = 	\underbrace{\begin{bmatrix}
        0 &  y_{t-1}
    \end{bmatrix}}_{H_t}\begin{bmatrix}
h_t \\
\beta_{1,t}
\end{bmatrix}  + \varepsilon_t, & \varepsilon_t & \sim {\mathcal N}(0,R), R = [\hat h_{t|t-1}]
\end{align*}
where the control input is defined by $u_{t-1} := [1, \hat e_{t-1}^2\Lambda_{(\hat e_{t-1})>0}, \hat e_{t-1}^2\Lambda_{(\hat e_{t-1})<0}]^\top$ and $\hat e_{t}$ is the KF residual, i.e. $\hat e_{t} := y_t - \hat \beta_{1,t|t-1} y_{t-1}$ at time $t$; see more detail in~\cite{1991:Hall:note}. The $\Lambda_{(\cdot)>0}$ takes the value $1$ when the last period's error is positive, and the value $0$ otherwise. Similarly $\Lambda_{(\cdot)<0}$ takes the value $1$ if the last period's error is negative, and the value $0$ otherwise. The unknown state vector $[h_t,\beta_{1,t}]^{\top}$ and system parameters $\theta = [\omega,a_1^{+},a_1^{-},b_1, \sigma_{w_1}^2]$ are to be estimated from the data, $Y_0^T:=\{y_t\}_1^T$ by using the methodology presented in~\ref{sec:appendix}.

Finally, the A-GARCH(1,1) in equations~\eqref{Asym2:yk1}~-- \eqref{Asym2:garch} has the following state-space form:
\begin{align*}
\underbrace{
\begin{bmatrix}
h_t \\
\beta_{1,t}
\end{bmatrix}
}_{x_t}  & =  \underbrace{
\begin{bmatrix}
b_1 & 0 \\
0 & 1
\end{bmatrix}}_{F}
\underbrace{
\begin{bmatrix}
h_{t-1} \\
\beta_{1,t-1}
\end{bmatrix}
}_{x_{t-1}}
 +
\underbrace{
	\begin{bmatrix}
		\omega  & a_1 & a_1^{+}\\
        0 & 0 & 0
	\end{bmatrix}}_{B}
\underbrace{
	\begin{bmatrix}
		1 \\
       \hat e_{t-1}^2\\
       (\hat e_{t-1}^{+})^2
	\end{bmatrix}}_{u_{t-1}} + 	\underbrace{\begin{bmatrix}
        0\\
        1
    \end{bmatrix}}_{G}
	\begin{bmatrix}
        w_{t}
    \end{bmatrix}, & w_t  & \sim {\mathcal N}(0,Q), Q = [\sigma_{w_1}^2], \\
    y_t  &  = 	\underbrace{\begin{bmatrix}
        0 &  y_{t-1}
    \end{bmatrix}}_{H_t}\begin{bmatrix}
h_t \\
\beta_{1,t}
\end{bmatrix}  + \varepsilon_t, & \varepsilon_t & \sim {\mathcal N}(0,R), R = [\hat h_{t|t-1}]
\end{align*}
where the control input is defined by $u_{t-1} := [1, \hat e_{t-1}^2,  (\hat e_{t-1}^{+})^2]^\top$ and $\hat e_{t}$ is the KF residual, i.e. $\hat e_{t} := y_t - \hat \beta_{1,t|t-1} y_{t-1}$ at time $t$; see more detail in~\cite{1991:Hall:note}. Additionally, $\hat e_{t}^{+} = \max \{\hat e_t, 0\}$. The unknown state vector $[h_t,\beta_{1,t}]^{\top}$ and system parameters $\theta = [\omega,a_1,a_1^{+},b_1, \sigma_{w_1}^2]$ are to be estimated from the data, $Y_0^T:=\{y_t\}_1^T$ by using the methodology presented in~\ref{sec:appendix}.

\end{document}